\mag=\magstephalf
\pageno=1
\input amstex
\documentstyle{amsppt}
\TagsOnRight
\interlinepenalty=1000
\NoRunningHeads

\pagewidth{16.5 truecm}
\pageheight{23.0 truecm}
\vcorrection{-1.0cm}
\hcorrection{-0.5cm}
\nologo

\NoBlackBoxes

\def\tvskip{\vskip 0.5 cm}

\define \ii{{\roman i}}

\define \bee{{\bold e}}
\define \bEE{{\bold E}}
\define \tr{{\roman {tr}}}

\define \dd{{\roman d}}

\define \RR{\Bbb R}
\define \EE{\Bbb E}
\define \ZZ{\Bbb Z}
\define \MM{\Cal M}
\define \AAA{\Cal A}
\define \EEE{{\Cal E}}

\define \NN{\Cal N}

\define \LL{\Cal L}
\define \III{{\Cal I}}

\define \CC{\Bbb C}

\define \conf{\roman{conf}}

\define \SS{\Cal{S}}
\define \fg{\frak g}
\define \tfg{\widetilde{\frak g}}
\define \tg{\widetilde{ g}}

\define \Hom{\roman{Hom}}
\define \SAnn{\Cal{A}nn}
\define \SHom{\Cal{H}om}
\define \SEnd{\Cal {E}nd}
\define \SGL{\Cal{G}\Cal{L}}

\define \GL{\roman{GL}}

\define \fvskip{\vskip 1.0 cm}
\define \CM{{\Cal C^\omega_M}}
\define \DM{{\Cal D_M}}
\define \DMod{{\Cal D}}
\define \CMinf{{\Cal C^\infty_M}}

\define \CEE{{\Cal C^\omega_{\EE^n}}}

\define \CTS{{\Cal C^\omega_{T_S}}}
\define \for{{\roman{for}}}
\define \CMC{{\Cal C^{\omega \Bbb C}_M}}

\define \CS{{\Cal C^\omega_S}}

\define \CT{{\Cal C^\omega_T}}
\define \RRM{{\Bbb R_M}}
\define \CCM{{\Bbb C_M}}
\define \gr{{\roman{gr}}}
\define \Lie{{\roman{Lie}}}

\define \Con{{\roman{Con}}}
\define \hor{{\roman{hor}}}

\define \Mod{{\roman{Mod}}}

\define \dalpha{{\dot \alpha}}
\define \dbeta{{\dot \beta}}
\define \dgamma{{\dot \gamma}}

\define \tnablaM{{\nabla^{\roman{SA}}_M}}
\define \tnablaT{{\nabla^{\roman{SA}}_{T_S}}}
\define \tnablaS{{\nabla^{\roman{SA}}_{S}}}
\define \XiM{{\Xi^{\roman{SA}}_M}}
\define \XiT{{\Xi^{\roman{SA}}_{T_S}}}


\topmatter
\title
\endtitle
{\centerline{\bf{ Submanifold Differential Operators  in $\Cal D$-Module 
Theory I :}}

{\centerline{\bf{Schr\"odinger Operators}}}

{}
{}
\vskip 0.5 cm
{\centerline{ Shigeki. Matsutani }

{\centerline{ 8-21-1 Higashi-Linkan Sagamihara 228-0811 JAPAN }

\address 
\endaddress
\abstract 

For this quarter of century,  differential operators in a lower dimensional 
submanifold embedded or immersed in real $n$-dimensional euclidean space 
$\EE^n$ have been studied as quantum mechanical models, which are realized as
restriction of the operators in $\EE^n$ to the submanifold. For this decade,
the Dirac operators in the submanifold have been investigated in such a scheme
,
which are identified with operators of the Frenet-Serret relation for a space
curve case and of the generalized Weierstrass relation for a conformal surface
case. These Dirac operators are concerned well in the differential geometry,
since they completely represent the submanifolds. In this and a future series
of articles, we will give mathematical construction of the differential 
operators on a submanifold in $\EE^n$ in terms of $\DMod$-module theory and 
rewrite  recent results of the Dirac operators mathematically. In this article
,
we will formulate Schr\"odinger operators in a low-dimensional submanifold in
$\EE^n$.

\centerline{{\bf MCS Codes:} 32C38, 34L40, 35Q40}

\centerline{{\bf Key Words:} Laplacian, Schr\"odinger operator, Submanifold,
$\Cal D$-Module}

\endabstract

\endtopmatter

\vskip 0.5 cm
{\centerline{\bf{\S 1. Introduction}}

Recently it becomes recognized that the Dirac operators play important roles
in  geometry {\it e.g.}, differential, algebraic, arithmetic geometry and so
on.  Pinkall  gave an invited talk in the international congress of 
mathematicians in 1998 on the relation between  immersed surfaces in three 
and/or four dimensional euclidean space $\EE^n$, ($n=3, 4$) and  Dirac 
operators, which was worked with Pedit [PP]. They constructed quaternion 
differential geometry and reduced the Dirac operators, which exhibit the 
geometrical properties of the surface. The Dirac operators of $\EE^n$  ($n=3,
4$) also had been discovered by Konopelchenko in studies on geometrical 
interpretation of soliton theory [Ko1, 2, KT] and by Burgess and Jensen [BJ]
and me [Mat3, Mat4] in the framework of the  quantum physics. Further on case
of $\EE^3$, Friedrich obtained it by investigation of spin bundle [Fr].

Our Dirac operator is purely constructed in analytic category as we will show
and is directly related to index theorems  [Mat2, Mat3, TM]. Thus I believe 
that it is important to reformulate our works in the framework of pure 
mathematics and  to translate them for mathematicians. In this and a future 
series of articles [II], I will mathematically  formulate the canonical 
Schr\"odinger operator and Dirac operator on a submanifold in $\EE^n$. Indeed,
there have appeared similar studies [DES, FH, RB] only on the Schr\"odinger 
operator case but it does not look enough to overcome several obstacle between
physics and mathematics. 

The submanifold quantum mechanics, which  I called, was opened by Jensen and
Koppe in 1971 [JK] and  rediscovered by da Costa in 1982 [dC]. They considered
a quantum particle confined in a subspace in our three dimensional euclidean
space $\EE^3$ which can be regarded as a low dimensional submanifold by taking
a certain limit. (Since confinement of quantum particle into a subspace is 
realized in a certain case [DWH], their investigation is not so fictitious.)
They found a canonical Laplacian by constructing the Schr\"odinger equation in
the submanifold, which differs from the ordinary Beltrami-Laplace operator [dC
,
JK]: For a surface embedded in $\EE^3$ case, the submanifold Laplacian 
$\Delta_{S \hookrightarrow \EE^3}$ is expressed by
$$
   -\Delta_{S \hookrightarrow \EE^3}:= -\Delta_S -(K-H),
$$ 
where $\Delta_S $ is a Beltrami-Laplace operator,
$K$ is Gauss curvature and $H$ is the mean curvature of $S$.
For a curve $C$ in $\EE^3$ case, we have 
$$
   -\Delta_{C \hookrightarrow \EE^3}:= -\Delta_C -\frac{1}{4}k^2,
$$ 
where $\Delta_C $ is a Beltrami-Laplace operator on $C$
and $k$ is curvature of the curve.

However submanifold quantum mechanics needs very subtle treatments. In fact,
there are several different types of theories of quantum mechanics for 
submanifolds. For example, it is well-known that restriction of quantum 
particle can be performed using Dirac constraint quantization scheme [Dir2].
Let an equation $f=0$ represent a hypersurface in $n$-dimensional $\EE^n$.
We can apply the Dirac constraint scheme [Dir2] to this system with a 
constraint $f=0$. Alternatively we can also deal with $\dot f=0$ constraint,
where dot means derivative in time $t$ [INTT]. These results differ; $\dot
f=0$ case agrees with the results of Jensen and Koppe [JK] and da Costa [dC]
whereas $f=0$ does not. Since the results of Jensen and Koppe [JK] and da 
Costa [dC] connect with fruitful results such as the generalized Weierstrass
equations as we will show in the introduction in [II], $\dot f=0$ case is 
very natural but $f=0$ case is not. In fact as a physical problem should be
determined by local information, the constraint $\dot f=0$ consists only of
local data whereas $f=0$ contains global information and is a fancy 
constraint from physical viewpoint [INTT, Mat1]. 

Accordingly when we make a theory of a submanifold quantum mechanics, 
we must pay many attentions on its treatment.

The self-adjoint operator is one of objects to need attentions. In order to
show importance of self-adjoint operator rather than canonical commutation
relation (or generating relation of Weyl algebra), let us consider a radial
differential operator of polar coordinate in $\EE^n-\{0\}$ [Dir1,SM]. For 
$n=2$ case, a point of $\EE^2-\{0\}$ is expressed as $\bold x=( r,\theta)$,
$\theta \in [0, 2\pi)$ and $r \in \RR_{>0}$ $:=\{r \in \RR |\ r>0\}$; $\RR$
is a set of real number. Then even though 
$\sqrt{-1}\partial_r:=\sqrt{-1}\partial/ \partial r$ is satisfied with the
canonical commutation relation $[\sqrt{-1}\partial_r , r ]:=
\sqrt{-1}\partial_r r -r \sqrt{-1}\partial_r =\sqrt{-1}$, the operator 
$\sqrt{-1}\partial_r$ is not self-adjoint if we use induced
 metric of $\EE^2$;
$$
    \int_0^\infty r d r\overline{\psi_1(r)} \sqrt{-1}\partial_r  
    \psi_2(r) \neq 
    \int_0^\infty r d r\overline{(\sqrt{-1}\partial_r \psi_1(r) ) 
    }\psi_2(r),
$$
where  $\psi_a$ $(a=1,2)$ which is a smooth function over $\RR_{>0}$ and  
vanishes at origin $r=0$ and asymptotically at $r=\infty$. (Here we have 
assumed that a function $\varphi$ over $\EE^2-\{0\}$ is expanded as $\varphi
=$ $\sum \psi^{(n)} \exp( \sqrt{-1} n \theta )$, where $\psi$'s are depend onl
y
on $r$.) We can hermitianize it by two ways [Dir1,dC]. First we will define an
operator $\sqrt{-1}\nabla_r := \sqrt{-1}(\partial_r + 1/(2 r))$, then 
$[\sqrt{-1}\nabla_r, r] = \sqrt{-1}$ and 
$$
\int_0^\infty r d r\overline{\psi_1(r)}\sqrt{-1} \nabla_r \psi_2(r) = 
\int_0^\infty r d r\overline{(\sqrt{-1}\nabla_r \psi_1(r) )} \psi_2(r).
$$
Second is that we will deform the function space such that 
$\Psi_a:= \sqrt{r}\psi_a$ and we obtain the relation,
$$
\int_0^\infty  d r\overline{\Psi_1(r)}\sqrt{-1} \partial_r  \Psi_2(r) = 
\int_0^\infty  d r\overline{(\sqrt{-1}\partial_r \Psi_1(r) )} \Psi_2(r).
$$
This transformation is well-known in  theory of ordinary differential equation
s
as a transformation from non-self-adjoint form to the self-adjoint one of 
Sturm-Liouville operator case (p.424 in [Arf].) It is not difficult to prove
that there is a ring isomorphism between a differential ring generated by 
$\sqrt{-1} \nabla_r$ and represented a vector space $\{\psi\}$ and that 
generated by $\sqrt{-1} \partial_r$ and represented over $\{\Psi\}$. In this
article we will consider them in the framework of $\DMod$-module.

Even though, $\sqrt{-1}\nabla_r$ is self-adjoint, $\sqrt{-1}\nabla_r$ is not
observable in general; the ordinary Schr\"odinger operator can not be expresse
d
by the spectral decomposition using eigen vectors of operator 
$\sqrt{-1}\nabla_r$ [Dir1]. (For example, radial momentum in hydrogen atom can
not observed.) However by an appropriate confine potential with infinite 
height, we can confine a quantum particle in a thin ring or thin surface; 
mathematically speaking, we can impose a Dirichlet boundary condition and 
restrict support of functions or the domain of the Schr\"odinger operator [dC,
DES, JK, SM]. Then normal mode for the codimensional direction in the ordinary
Schr\"odinger equation, which is expressed by (bilinear of) $\sqrt{-1}\nabla_r
$
is well-defined and $\sqrt{-1}\nabla_r$ behaves as momentum operator. (In 
meso-scopic quantum mechanical system, normal mode can be observed as subband
state [DWH].) 

Let us take a squeezing limit so that thickness of subspace or support of 
functions is negligible and the subspace can be regarded as a lower dimensiona
l
submanifold $S$ itself. Then we can argue a system of differential operators
defined over the subspace; a self-adjoint differential operator along the 
normal direction can be constructed similar to $\sqrt{-1}\nabla_r$. By 
integrating the hamiltonian over the normal mode, we obtain a submanifold 
quantum mechanics. This procedure is similar to the techniques to get the Thom
isomorphism using Berezin integration method [BGV, Y].

However it is not easy to justify such squeezing  limit including Dirac 
operators  using concept of Hilbert space even though Duclos, Exner and
{\u S}{\u t}ov\'i{\u c}ek attempted for Schr\"odinger operators [DES]: In 
squeezing limit,  we must evaluate divergence of eigenvalue of normal 
direction. 

Thus in this article, I will make an attempt to reformulate it using 
$\DMod$-module theory. In fact, Sato said that to study noncommutative system
sometimes needs an appropriate topology instead of ordinary topologies which
are used in the operator algebra, {\it e.g.}, weak topology in Hilbert space
[S]. Further in the submanifold quantum mechanics, we need a restriction of
the differential operator while the restriction is the most natural concept
in the sheaf theory and $\DMod$-module theory is based upon the sheaf theory.
Thus I believe that my attempt is more natural than others approaches to the
submanifold quantum mechanics. 

The $\DMod$-module theory was began by Sato as a algebraic analysis [Bjo, K].
For a noncommutative ring $\DMod_{\EE^n}$ of differential operators over 
$n$-dimensional parameter spaces and a given differential equation $P u=0$ for
$P\in \DMod_{\EE^n}$, let us consider a left $\DMod_{\EE^n} $-module
$$
    \Cal M_{\EE^n} = \DMod_{\EE^n}/\DMod_{\EE^n} P .
$$
Then the ring homomorphism of $\Cal M_{\EE^n}$ to a  function space such as 
analytic function space $C^\omega_{\EE^n}$ is ring isomorphic to solution spac
e
of the equation $P u=0$ [Bjo, Cou, HT]. For a more general case, above quotien
t
space is replaced with a coherent module. In the $\DMod$-module theory, the 
differential operators on a submanifold have been studied in detail. However
in these studies, our differential operators in submanifold quantum mechanics
have not ever appeared as long as I know.

The algorithm to construct the submanifold quantum mechanics
in $\DMod$-module theory for a hypersurface $S$ in $n$-dimensional
euclidean space $\EE^n$ is as follows.

\roster

\item We construct a quantum equation and its related $\DMod_{\EE^n}$-module
$\Cal M_{\EE^n}$ in $\EE^n$ with the natural metric, {\it e.g.}, for the free
Schr\"odinger equation, $-\Delta_{\EE^n}\psi=0$, $\Cal M_{\EE^n}=\Cal 
D_{\EE^n}^\CC/(\Cal D_{\EE^n}^\CC (-\Delta_{\EE^n}))$; Here $\Cal 
D_{\EE^n}^\CC:= \Cal C_{\EE^n}^\CC[\partial_1, \cdots, \partial_n]$, $\Cal 
C_{\EE^n}^\CC$ is complex valued analytic functions over $\EE^n$ and 
$\Delta_{\EE^n}$ is the Beltrami-Laplace operator in $\EE^n$.

\item We embed (or immerse) a real analytic hypersurface $S$ in  $\EE^n$, whic
h
is given by an equation $f=0$.

\item We find a local system along  a tubular neighborhood $T_S$ of the 
submanifold $S$ and calculate the inverse image $\Cal M_{T_S}$ of $\Cal 
M_{\EE^n}$ to  $T_S$ and $\Delta_{T_S}:= \Delta_{\EE^n} |_{T_S}$ for $T_S  
\hookrightarrow \EE^n$.

\item We find a self-adjoint operator $\sqrt{-1}\tnablaS^\perp$ along 
the normal direction of $S$ over an open set $U$ of $S$.

\item We define a quotient module $\Cal M_{S\hookrightarrow T_S}$ by
an exact sequence, 
$$
    0 \longrightarrow \Cal M_{T_S}|_S
     \overset{\sqrt{-1}\tnablaS^\perp}\to{\longrightarrow}
    \Cal M_{T_S}|_SS\longrightarrow\Cal M_{S\hookrightarrow \EE^n}
    \longrightarrow 0.
$$

\item $\Cal M_{S\hookrightarrow \EE^n}$ is a coherent $\DMod_S$-module related
to submanifold quantum mechanics over the submanifold $S$. We will define a 
submanifold quantum mechanical operator $\Delta_{S\hookrightarrow \EE^n}$ by
the exact sequence,
$$
  0 \longrightarrow \DMod_S 
  \overset{\Delta_{S\hookrightarrow T_S}}\to{\longrightarrow} \Cal D_S
   \to \Cal M_{S\hookrightarrow \EE^n} \longrightarrow 0.
 $$

\endroster

\vskip 0.5 cm

This scheme does not need any limit procedure and avoid the disease of 
divergence. (5) and (6) can be rewritten as follows if you choose a local 
coordinate system,

\roster
\item"(5')" We will define $\Delta_{S\to T_S}:= \Delta_{T_S}|_S$ and 
$\Delta_{S\hookrightarrow \EE^n} := \Delta_{S \to 
T_S}|_{\sqrt{-1}\tnablaT^\perp=0}$ for a standard form of $\Delta_{S\to T_S}$,
where all $\tnablaS^\perp$ are put right side of each terms of $\Delta_{S \to
T_S}$.

\endroster

Physically speaking, these processes are naturally performed when we introduce
the confinement potential along the submanifold with the same thin thickness
[dC, JK].

Contents are as follows. In \S 2, we will quickly review $\DMod$-module theory
and sheaf theory. In \S 3, we will define the adjoint operator using Hodge $*$
product, though it can be defined using extension in the cohomology theory. We
will introduce the half-from and the self-adjoint momentum  operator.

In \S 4, we will define the Schr\"odinger  operator in a
lower dimensional submanifold in $\EE^n$ and give theorems.

\vskip 1.0 cm

\centerline{\bf Acknowledgment}
\tvskip 
It is acknowledged that Prof.~K.~Tamano have taught me algebraic topology, 
sheaf theory and differential geometry   for this decade and Prof. Y.~\^Onishi
have discussed on the number theory, algebraic geometry and been consulted my
mathematical approach for this decade. 

I am grateful to Prof.~A.~Kholodenko for  so many encouragements and 
discussions by using e-mails, to Prof.~B.~L.~Konopelchenko for kind letters to
encourage these works and to Prof. F.~Pedit for his interesting in this work. 

I thank Dr.H.~Tsuru, Prof.~S.~Saito, Prof.T.~Tokihiro, Prof.~M.~Nakahara, 
Prof.A.~Suzuki, Prof.S.~Takagi, Prof.~H.~Kuratsuji, Prof.~K.~Sogo, W.~Kawase
and H.~Mitsuhashi for helpful discussions and comments in early stage of this
study.

I thank to Prof.~Y.~Ohnita, Prof.~M.~Guest, Dr.~R.~Aiyama and 
Prof.~K.~Akutagawa for inviting me their seminars and for critical discussions
.


\vskip 1.0 cm
{\centerline{\bf{\S 2.  Foliation and $\DM$-Module }}

\vskip 0.5 cm

Although we will not construct a theory in the category of differential 
geometry whose morphism is ($\Cal C^\infty$)-diffeomerphism, we are concerned
with differential geometry and physical system rather than merely algebraic 
structure. Accordingly we will treat {\it real}  analytic objects in this 
article. 
 
Let $M$ be a real $n$-dimensional analytic manifold without singularity and 
$\CM$ is a structure sheaf of real analytic functions over $M$; ($M,\CM)$ is
an algebraized analytic manifold. A sheaf $\Cal S$ (of sets) over a topologica
l
space $X$ is characterized by a triple $(\Cal S, \pi, X)$ because it is define
d
so that there is a local homeomorphism $\pi: \Cal S \to X$ [Mal1]. In this 
article, we choose $M$ as such a topological space $X$. Economy of notations
make us to denote $\Cal S:=(\Cal S, \pi, M)$  for abbreviation. Further we wil
l
write a set of (local) sections of $\Cal S$ over an open set $U$ of $M$ as 
$\Cal S(U)$ or $\Gamma(\Cal S, U)$. Using the category equivalence between 
category of sheaves and category of complete presheaves (due to theorem 13.1
in [Mal1]), we will mix them here.

\subheading{Notations 2.1 (Sheaves) }[Mal1, Mal2, Bjo] 

{\it \roster 

\item Complexfication of $\CM$ is denoted by $\CMC$.

\item Let us denote a sheaf $\CC_M$
consisting of a locally constant functions $\CC$ over $M$.
Similarly $\RR_M$ $\ZZ_M$,  ${\ZZ_2}_M$ and 
so on. 

\item $1_M$ means a unit element
sheaf of multiplicative group sheaf $\CC_M^\times$ of $\CC_M$.

\item If  $\Cal E$ is a finite rank of locally free $\CM$-module,
we call it a vector sheaf.

\item The tangent sheaf over $M$,  which is a vector sheaf, is written as
$$
 \Theta_M :=Der_{\RRM}(\CM):=\{\theta \ | \ 
 \theta(f g ) = \theta(f) g + f \theta(g),  \ 
 (f,g \in \Gamma(U,\CM)), \ U \subset M (\text{open}) \}.
$$
{\rm{(p.18 1.2.7 in [Bjo])}}

\item The complex valued tangent sheaf is
$$
\Theta_M^\CC =Der_{\CCM}(\CMC):=\Theta_M \otimes \CC.
$$ 
{\rm{(p.18 1.2.7 and p.281 in [Bjo])}}

\item Let us denote the sheafification of homomorphism of $\Cal A$-module 
sheaves $\MM$ and $\MM'$  over $M$ by $\SHom_{\Cal A}(\MM,\MM')$. Similarly th
e
sheafification of endomorphism of  $\Cal A$-module sheaves $\MM$  over $M$ by
$\SEnd_{\Cal A}(\MM)$. {\rm{(p.133-p.144 in [Mal1])}}

\item A locally constant sheaf $\LL_M$ over $M$, called local system, is 
defined so that its stalks are finite 
dimensional real vector spaces $\RR^m$. {\rm{(p.22 in [Bjo])}}.

\item A group sheaf generated by the presheaf,
$$
    U \to \GL(n,\CM)(U),
$$
is denoted by general linear group sheaf $\SGL(n,\CM)$, where 
$U$ is an open set of $M$ and $\GL(n,\CM)$
is general linear group for $\CM(U)$-valued $n$-matrix.
 {\rm{(p. 285 in [Mal1])}}

\item 
A sheaf homomorphism as group from $\SGL(n,\CM)$ to $\CC_M^{\times}$ is
denoted by $\det$. {\rm{(p.294-p.295 in [Mal1])}}

\endroster
}

\tvskip

\subheading{Proposition 2.2}

{\it For each point $p \in M$, there exists an open neighborhood $U_{p}$ aroun
d
$m$ with a local coordinate system $\{x_i,\partial_i\}_{1\leq i \leq n}$ 
satisfied with,
$$
    x_i \in \Gamma(U, \CM ),\quad
    \Theta_M(U) = \bigoplus_{i =1}^n \Gamma(U, \CM 
    )\partial_i ,
    \quad[\partial_i,x_j] \equiv \partial_i x_j- 
    x_j\partial_i= \delta_{i j } .
$$ 
}

\demo{Proof} see [Bjo] p.11 proposition 1.1.18 and p.17 
remark under definition 1.2.2. \qed
\enddemo

\tvskip
\subheading{Definition 2.3} (1.2.8 in [Bjo])

{\it In a chart, we write sections $\delta$ and $\delta'$ in $\Theta_M(U)$, $
\delta= \sum_i f_i(x) \partial_i$, $\delta'= \sum_i g_i(x) \partial_i$, where
$U$ is an open set of $M$, $f$'s and $g$'s are in $\CM(U)$. The commutator is
defined by
$$
    [ \delta, \delta'] = \sum_{u,v} 
    \left(    f_v \partial_v (g_u) \partial_u 
    -g_u \partial_u (f_v) \partial_v \right) ,\quad
    [\delta, g] = \sum_{u,v} 
    f_v \partial_v g.
$$
}

\subheading{Definition 2.4 (Differential Ring Sheaf)}[Bjo, TH]

{\it 
We will denote the subring sheaf $\DM$ of $\SEnd_\RRM(\CM)$, which is generate
d
by the complete presheaves $\Gamma(\CM,U)$ and $\Gamma(\Theta_M,U)$ and has 
local expression,
$$
           \Gamma(U, \DM)  = \bigoplus_{\alpha \in \Bbb N^n}
            \Gamma(U, \CM )\partial^\alpha ,
$$
where $ \partial^\alpha = \partial_1^{\alpha_1} 
 \partial_2^{\alpha_2} \cdots \partial_n^{\alpha_n}\in\Gamma(\Theta_M,U)$.
 We will define $\DM^\CC$ as $\DM \otimes \CC$.
}

\subheading{Remark 2.5 (Standard From Representation)}[Bjo]

The local expression is based upon the {\it standard form representation},
which is, for any $P \in  \Gamma(U, \DM)$, given as, 
$$
           P = \sum_{\alpha \in \Bbb N^n} a_\alpha\partial^\alpha ,
           \quad a_{\alpha} \in \CM(U).
$$

\subheading{Definition 2.6 (Filter)} [Bjo, TH]

{\it 
Let us define the filter $F$ associated with rank as a sheaf morphism by 
complete presheaf generated by the sections for an open set $U$ of $M$,
$$
    (F_l \DM)(U) = \sum_{|\alpha|\leq l, \ \alpha \in \Bbb N^n} 
     \Gamma(U, \CM ) \partial^\alpha,
$$
where $|\alpha| = \sum_i \alpha_i$.
Further for an open set $V$ of $M$, the filter $F_l$ of
$\DM$ is
defined as
$$
    (F_l \DM)(V):=
        \{ P \in \DM(V)\ | \ 
    \rho_{U V} P \in F_l(\DM(U)), \ (U \subset V)\},
$$
where $\rho_{U V}$ is a restriction associated with
the complete presheaf $\DM$.}

\subheading{Proposition 2.7}

{\it
\roster
\item  $F_l$ is increasing filter  and defined over $M$;
$$
    \DM = \bigcup_{m \in \Bbb N} F_l \DM,\quad
    F_{l}\DM\subset F_{m} \DM , \ \roman{for} \ l <m. 
$$

\item 
$F_0 \DM = \CM,$ $(F_l \DM ) (F_m \DM) = F_{l+m} \DM$ as 
$\CM$-module.

\item
For $P \in F_l \DM $ and $Q\in F_m \DM$, 
$[P,Q]\in F_{l+m-1} \DM $,
where $F_l \DM =0 $ for $l<0$; symbolically
$$
    [(F_l \DM ), (F_m \DM)] \subset F_{l+m-1} \DM.
$$

\item
We will define a symbol, 
$\gr \DM := \bigoplus_{l=0} (F_l \DM)/(F_{l-1} \DM)$.
Then $\gr \DM$ is a commutative ring and is regarded
as a commutative $C^*$ algebra. 

\endroster
}

\demo{Proof} [Bjo] p.17-p.18 and [TH] \qed
\enddemo

\tvskip

\tvskip
\subheading{Definition 2.8 (Differential Form)} (p.240 in [Mal1])

{ \it
\roster

\item   We will define the duality 
$ \Omega_M^1 := \SHom( \Theta_M, \CM)$ generated by
$\omega(\theta) \in \CM(U)$ for 
$\omega \in \Gamma(U,\Omega_M^1)$ and $\theta
\in \Gamma(U,\Theta_M)$ where $U$ is an open set of $M$.

\item 
Let $\Omega_M$ be a sheaf of differential form as  
a graded commutative ring generated by  $ \Omega_M^1$,
$$
 \Omega_M= \bigcup_{m \in \ZZ_{\ge0}}\Omega_M^p,
$$
where $ \Omega_M^p$ is set of $p$-forms and $\Omega_M^q=0$
$(q>n)$ vanishes: $\ZZ_{\ge0}:=\{n \in \ZZ\ | \ n \ge 0\}$

\item
The exterior derivative  is expressed by 
$ d: \Omega_M^p \to \Omega_M^{p+1}$. 

\endroster
}

\tvskip
\subheading{Definition 2.9 (Lie Derivative)}[TH,W]

{ \it 
For  an open set $U$ of $M$ and $\theta \in \Theta_M(U)$, 
the Lie derivative $\Lie_\theta$ 
is defined as follows:

\roster
\item for $f\in \Gamma(U,\CM)$, $
     \Lie_{\theta_1}(f) = \theta_1f. $
     
\item For  $\theta_a \in \Theta_M(U)$ ($ a =1,2$), $
    \Lie_{\theta_1}(\theta_2) = [\theta_1,\theta_2].$
    
\item For $\theta_a \in \Theta_M(U)$ ($ a =1,2$) and
 $\omega \in \Omega_M^1(U)$, $
\Lie_{\theta_1}(\omega(\theta_2)) =
    (\Lie_{\theta_1}(\omega))(\theta_2) - 
    \omega(\Lie_{\theta_1}(\theta_2))
    =\theta_1(\omega(\theta_2)) - 
    \omega([\theta_1,\theta_2])$.
    
\item For $\theta_a \in \Theta_M(U)$ ($ a =1,\cdots,n$) 
and $\omega \in \Omega_M^n(U) $, $
\Lie_\theta(\omega(\theta_1,\theta_2\cdots,\theta_n))    
        $ $=$ $\theta(\omega(\theta_1,\theta_2\cdots,$ $\theta_n)) 
      $ $-\sum_{i=1}^n 
       \omega(\theta_1,\cdots,[\theta,\theta_i],\cdots,\theta_n).$
\endroster
}

\subheading{Proposition 2.10}
{\it
\roster

For  an open set $U$ of $M$, $f \in \CM(U)$
 and $\omega\in \Omega^n_M(U)$,
$
    \Lie_{f \theta} \omega = \Lie_\theta f\omega .
$

\endroster
}

\demo{Proof} Direct computation shows it [TH]. \qed\enddemo

\subheading{Definition 2.11 (Right Action)} ([TH], [Bjo])

{ \it

The left action of $\theta$ of section of 
$\Theta(U)$ for an open sets $U$ of $M$
is defined by  
$$
    \omega\cdot \theta:= - (\Lie_{\theta}) \omega,
$$
where $\omega \in \Omega^n_M(U)$.
}

\subheading{Definition 2.12 (Integrable Connection)}( 1.2.10 in [Bjo])
{\it

Let $\Cal F$ be $\CM$-Module. Put
$
    \III (\Cal F):= \SHom_{\RR_M}( \Cal F, \Cal F ).
$ 
A global section $\nabla$  of  
$\SHom_{\CM}( \Theta_X, \III(\Cal F) )$
is called an integrable connection on $\Cal F$ if the following
holds for any open set $U$ of $M$:
\roster

\item For $\alpha \in \CM(U)$ and $f \in \Cal F(U)$, 
$
    \nabla_\delta( \alpha f)
     = \delta \alpha \cdot f + \alpha \nabla_\delta (f).
$

\item
$
\nabla_{[\delta,\delta']} = [\nabla_\delta,\nabla_{\delta'}].
$

\item $\nabla$ is involutive, i.e., $ [\nabla_\delta,\nabla_{\delta'}]$
 is a commutator in $\III(\Cal F)$.
\endroster}

\subheading{Definition 2.13 (Category)} (1.1.24 and 1.2.11  in [Bjo])
{\it
\roster 

\item Denote by $\Mod^L(\DM)$ the category of left $\DM$-module.

\item Denote by $\Mod(\CM)$ the category of  $\CM$-module.

\item $\for$ is the forgetful functor from $\Mod^L(\DM)$ to $\Mod(\CM)$.

}

\tvskip

\subheading{Proposition 2.14 }(1.2.11  in [Bjo])
{\it 

Consider the category $\AAA$ whose objects consist of pair $(\NN, \nabla)$,
where $\NN \in \Mod(\CM)$ and $\nabla$ is an integral connection on $\NN$.
Morphisms are defined as follows:
$$
\Hom_{\AAA} ( (\NN, \nabla), (\NN', \nabla') )
=\{ \ \varphi\in \SHom_{\CM}(\NN,\NN') \ | \ 
\nabla' \circ \varphi = \varphi \circ \nabla \ \}.
$$

Then $\Mod^L(\DM)$ and $\AAA$ are category equivalent by the functor
$\mu: \Mod^L(\DM) \to \AAA$ for which $\mu(\MM) = (\for(\MM), \nabla)$ where
$\nabla_\delta (m ) = \delta(m)$ for every $\delta \in \Theta_M$ and
$m \in \MM \in \Mod^L(\DM)$.

}

\demo{Proof} [Bjo] p.19  Theorem 1.2.12.\qed \enddemo

This means that we can find an object in $\Mod^L(\DM)$ for any an integrable
connection set $(\nabla, \NN)$.

\tvskip

\subheading{Definition 2.15 (Horizontal Section)} (1.3.7 in [Bjo])

{\it   
The left annihilator of $1_M$ be denoted by.  
$$
     \SAnn_{\DM}(1_M):= \{ P \in \DM \ | \ P( 1_M ) = 0 \ \}.
$$
For every $\MM \in \Mod^L( \DM ) $, we introduce $\RR_M$-module,
called the sheaf of horizontal sections of $\MM$,
$$
    \hor( \MM ):=\SHom_{\DM} ( \DM/\SAnn_{\DM}(1_M), \MM) .
$$
Then we have
$$
    \MM \approx \CM \otimes_{\RR_M} \hor(\MM), \quad
    \DM/\SAnn_{\DM}(1_M) \approx F_0 \DM\approx \CM.
$$
} 

\tvskip
\subheading{Proposition 2.16 (Connections)} (1.3.9 in [Bjo])

{\it  
\roster

\item For a local system $\LL_M$, the left $\DM$-module $\CM \otimes_{\RR_M}
\LL_M$ is denote by $\Con(\LL_M)$. Let 
a set of local systems $\LL_M$ be $\frak L_M$
and  $\Con(\DM)$ $:=\{\Con(\LL_M)\ | \LL_M\in\frak L_M \}$. 
$\frak L_M$ and $\Con(\DM)$ are category equivalent due to the 
relation; for $\LL_M$, $\LL_M'$ $\in\frak L_M$ 
$$
    \SHom_{\DM}( \Con(\LL_M), \Con(\LL_M')) 
    =\SHom_{\RR_M}(\LL_M, \LL_M')  .
$$

\item For a local system $\LL_M$, there exists a correspondence,
$$
    \SHom_{\DM}( \DM/ \SAnn_{\DM}( \LL_M ), \CM) = \LL_M,
$$
where
$$
     \SAnn_{\DM}(\LL_M):= \{ P \in \DM \ | \ P( \LL_M ) = 0 \ \}.
$$
\endroster
} 

\demo{Proof} [Bjo] p.22-p.23.\qed \enddemo 

If $\LL_M$ is a subset of $\Theta_M$, this theorem is essentially the same as
the Frobenius integrable theorem [Mal2]. The local system can be regarded as
distribution in the terminology of differential geometry [W]. \fvskip

Let $S$ be a real analytic submanifold without singularity of the manifold $M$
.
For a point $s \in S$, there is an open neighborhood of $s$, there exist  real
analytic functions $Q_1,\cdots,Q_d \in \CM$ such that
$$
    S\cap U =\{ s\subset U \ |\ Q_1(s)=\cdots=Q_d(s) = 0 \},
$$
where $d=n-k$.

Let the natural embedding be expressed by $\iota_S : S \hookrightarrow M$.
Then for a sheaf $\Cal F$ over $M$, we will define
$\Cal F|_S := \iota_S^{-1} \Cal F$.

\fvskip

\subheading{Definition 2.17 (Inverse Image of $\DM$-module)}

{\it
The sheaf $\DMod_{S \to M}$ over $S$ is defined as 
$\DMod_{S \to M} = \CS \otimes _{\CM|_S} (\DM|_S)$.
}
\fvskip

If we will assign the local coordinate of open neighborhood
of a point $s \in S\subset M$, $q^1=q^2=\cdots=q^d=0$.
Locally $\DMod_{S \to M}$ is exxpressed as,
$$
    \DMod_{S \to M}= \DM|_S/(q^1( \DM|_S)
    +q^2( \DM|_S)+\cdots+q^d( \DM|_S)).
$$

\fvskip

Next we will introduce the Riemannian metric in sheaf theory
along the line of the arguments of Mallios [Mal1,2].

\subheading{Definition 2.18 (Ordered Algebraized Space)}(Definition 8.1 in
[Mal1])

{\it
Let $\AAA$ be a real ring sheaf whose local sections are real valued. Suppose
the $(M,\AAA)$ is algebraized space and $\AAA^+$ is a subsheaf of the real rin
g
sheaf $\AAA$, whose local sections are positive real valued. $(M,\AAA, \AAA^+)
$
is an ordered algebraized space, viz, for a any local section $\lambda \in 
\RR^{>0}_M(U)$, $\lambda\AAA^+(U)  \subset \AAA^+(U)$, $\AAA^+(U) + \AAA^+(U)
\subset \AAA^+(U)$ and $\AAA^+(U)  \AAA^+(U) \subset \AAA^+(U)$.
}

\fvskip

Then it is obvious that there exists  an ordered algebraized space 
 $(M,\CM, \CM^+)$.

\subheading{Definition 2.19 (Inner product module)} (p.318 in [Mal1]) 
{\it

Suppose that $(M,\AAA, \AAA^+)$ is an order algebraized space 
and  $\EEE$ is $\AAA$-module.
We say that $(\EEE,\fg_\EEE)$ has an $\AAA$-valued inner product
and   $(\EEE,\fg_\EEE)$ is an inner product $\AAA$-module on $M$,
if a sheaf morphism $\fg_\EEE: \EEE\oplus\EEE \to \AAA$
 is satisfied following conditions:

\roster

\item $\fg_\EEE$  is a $\AAA$ bilinear morphism.

\item  $\fg_\EEE$ is positive definite; for any local section $s\in$
$\EEE(U)$ over an open set $U$, $\fg_\EEE(s,s)\in \AAA^+(U)$
such that $\fg_\EEE(s, s) =0$ if and only if $s=0$ in $\EEE(U)$.
 
\item   $\fg_\EEE$ is symmetric; for any local two sections $s$ and $t$
of $\EEE(U)$ over an open set $U$, $\fg_\EEE(s,t)=\fg_\EEE(t, s)$.

\endroster

}

\subheading{Definition 2.20 (Riemannian metric)} (p.320 in [Mal1]) 
{\it

We say that an order algebraized space $(M,\CM, \CM^+)$
 has a Riemannian metric $(\Theta_M,\fg_M)$ as 
  a sheaf morphism $\fg_M: \Theta_M \to \Omega_M^1$ if 
it is satisfied with following conditions

\roster

\item The sheaf morphism $\fg_M: \Theta_M \to \Omega_M^1$ is 
$\CM$-isomorphism.

\item Using the natural duality $\Theta_M \oplus \Omega_M^1 \to \CM$,
$\fg_M$ is extended to $\CM$-valued inner product 
$\fg_M: \Theta_M \oplus \Theta_M \to \CM$.

\item The action of $\theta \in  \Theta_M$ is defined by
$$
    \theta \fg_M( \theta_1,\theta_2) =
     \fg_M( \theta \theta_1,\theta_2) 
     + \fg_M( \theta_1,\theta \theta_2).
$$
\endroster}

Here we will note {\it {strictly fine sheaf} }, which is used in the definitio
n
of the Riemannian module in [Mal1] in the category for sheaves of smooth 
functions $\CMinf$. By using forgetful functor from the category of $\CM$ to
that of $\CMinf$, it turns out that image of the functor is a subset of 
strictly fine sheaves [Mal1]. Hence above definition is not contradict with 
that in [Mal1].

Further we will define a morphism 
$$
    \tfg_M : \Omega_M^1 \to \Theta_M,
$$
as $\tfg_M \circ \fg_M = id_{\Theta_M}$ and $\fg_M \circ \tfg_M = 
id_{\Omega_M^1}$. Then due to their duality, $\tfg_M: \Omega_M^1\oplus 
\Omega_M^1 \to \CM$.

Since analytic  manifold $M$ is given by a solution of a certain differential
equations and we assume that our considered manifold $M$ has no singular, from
proposition 2.16, it has a local coordinate system. Using local frame 
$(x^a)_{a=1,\cdots,n}$, as the duality is expressed by $< d x^a, \partial_s >
= \delta^a_b$ and $\Theta_x$ is decomposed as $\Theta_x = \sum_i C_x^\omega 
\partial_i$ where germs of $\Theta_M$, $C_x^\omega$ at $x$, we can express it
$$
    g_{M i,j}\equiv g_M(i,j):=\fg_M(\partial_i, \partial_j),
$$
and thus
$$
    \fg_M \equiv g_{M i,j} dx^i \otimes_{\CM} dx^j.
$$
Since $\fg_M$ can be realized as a section of $\SGL(n,\CM)$, there is a map 
$\det: \SGL(n,\CM) \to \CM$, we will denote $g_M := \det( g_{M}( i,j))$.

\subheading{Definition 2.21 (Normal Sheaf)}

{\it
Let us define the normal sheaf $\Theta_{S}^\perp$ as $\CS$-module
sheaf by the exact sequence,
$$
 0 \to \Theta_S \to \iota_S^{-1} \Theta_M \to \Theta_{S^\perp} \to 0,
$$
and the integrable connection of $\Theta_{S^\perp}$ as
$ \nabla_S^\perp \in \SHom_{\CS}( \Theta_S, \III(\Theta_{S^\perp}))$.
}

\subheading{Lemma 2.22 }
{\it

Let $\Omega^{1}_{S^\perp} := \SAnn_{\iota_S^{-1}\Omega_M^1}( \Theta_S)$
$=\{ \omega \in \iota_S^{-1}\Omega_M^1 \ | \ \text{for}
 \forall \theta \in  \Theta_S, \ \omega(\theta)=0\ \}$.
Then we have a natural correspondence
$$
	\tfg_M( \Omega^{1}_{S^\perp}) = \Theta_{S^\perp}.
$$

 By applying the proposition 2.14
and 2.16, we obtain the theorem.

\subheading{Theorem 2.23 (Existence of Tubular neighborhood)}
{\it

\roster
Suppose $(S, \CS)$ is a $k$-dimensional non-singular real analytic algebraized
submanifold of ($M,\CM)$: $\iota_S:S \hookrightarrow M$. We have a tubular 
neighborhood $(T_S,  \LL_{T_S}^\parallel)$ of $S$ satisfied with
following conditions,

\item $T_S$ is an open set of $M$, whose dimension is the same as $M$ as a
manifold, such that there are natural real analytic inclusions,
 $i_{S} :S\hookrightarrow T_S$, and 
$i_{T_S} :T_S\hookrightarrow M$, where 
$\iota_S\equiv i_{T_S}\circ i_{S}$.

\item There is a real analytic projection from $T_S$ to $S$, $\pi_{T_S}: T_S
\to S$ such that $\pi_{T_S} \circ i_S = id_S$.

\item The tangent sheaf of $\Theta_{T_S}$ is $\Theta_{T_S}=i_{T_S}^{-1} 
\Theta_M$ and has a direct decomposition as a $\CTS$-module,
$$
    \Theta_{T_S} = \Theta_{T_S}^\parallel\oplus \Theta_{T_S}^\perp,
$$
where $i_S^{-1}\Theta_{T_S}^\parallel =\Theta_{S}^\parallel$ and $
i_S^{-1}\Theta_{T_S}^\perp =\Theta_{S^\perp}$.

\item $\LL_{T_S}^\parallel$ is a local system over $T_S$ such that $\DM$-modul
e
$\MM_{T_S}^\parallel :=\SAnn_{\DMod_{T_S}}( \LL_{T_S}^\parallel )$ is generate
d
by $\Theta_{T_S}^\parallel$ and, whose local sections are given as 
$\Gamma(U,\LL_{T_S}^\parallel) = \Gamma(U,\RR_{T_S}^{n-k})$ for any open set
$U$ in $T_S$.

\endroster
}

\demo{Proof}
Using proposition 2.16 and the fact that $M$ has no singular, we can prove
them. However they are also proved in concepts in Frobenius integrability 
theorem [AM, Mal2, W]  and complete parallelism [Mal2] p.136. Later is
familiar for differential geometrists. \qed
\enddemo

\fvskip
\subheading{Remark 2.24 (Tubular neighborhood)} 

\roster

\item $ T_S$ has a 
foliation structure and $( \LL_{T_S}^\parallel, \pi_{T_S}, S)$ is a sheaf 
over $S$.

\item For the metric $\fg_M$ of $M$, $T_S$ and $S$ have
 the induced metric $i^{-1}_{T_S}\fg$ and $i_{S}^{-1}\fg$
respectively.

\item
We can define a local coordinate system of a open set $U$ of
$T_S$: $p \in U$, p is expressed by
$(s^1,\cdots,s^k,q^{k+1},\cdots,q^n ) $,
where $(s^1,\cdots,s^k)$ is a local coordinate system of $S$
by $\pi_{T_S} U$.
Then $\LL_{T_S}^\parallel(U)$ is characterized by a constant 
section $(q_{k+a})_{a=1,\cdots,n-k}$  $\in \RR^{n-k}_M(U)$ 
and exhibits a leaf of foliation.

\endroster

\vskip 1.0 cm
{\centerline{\bf{\S 3. Self Adjoint Operator}}
\fvskip

Let $\Mod^L \DM$ and $\Mod^R \DM$ be abelian category
of left and right $\DM$-modules.

\subheading{Proposition 3.1}

{\it
\roster
\item
For $\Cal M\in  \Mod^L \DM$,
if we define the action of tangent sheaf $\Theta_M$ to
$\Cal M \otimes_{\CM}\Omega_M^n$ as
$$
( m \otimes \omega) \theta:= -(\theta m) \otimes \omega
                          + m \otimes \omega \theta,
$$
where $m \otimes \omega \in \Cal M \otimes_{\CM}\Omega_M^n$, $\theta \in 
\Theta_M$, then $\Cal M \otimes_{\CM}\Omega_M^n$ can be regarded as the right
$\DM$-module.

\item $\Omega_M^n$ can be regarded as the isomorphism
as abelian category, which maps
$$
    \Omega_M^n:\Mod^L \DM\to \Mod^R \DM .
$$
\endroster}

According to the arguments of Mallios [Mal2, p. 343], we will
introduce the Hodge *-operator, volume element and
Beltrami-Laplace operator.
Hereafter we will sometimes use Einstein convention; we will
sum over an index if it appear twice in a term.

\subheading{Definition 3.2} [AM, Mal1, N, W]

{\it
Let $(M,\CM,\CM^+)$ be a non-singular analytic manifold 
endowed with the a Riemannian $\CM$-module $(\Theta_M,\fg_M)$. 
We have the Hodge operator as an element of 
automorphisim of $\Omega_M$,
$$
    *: \Omega_M^p \to \Omega_M^{n-p} .
$$
By local chart $U$, the metric is expressed by 
$\fg_M=g_{M i,j} d x^i \otimes d x^j$ over $U$,
there $\tg_M^{ij}$ is the inverse matrix of $g_{M i,j}$,
and its determinant is expressed by $g_M:= \det g_{M  i,j}$,
then the Hodge operator is represented by,
$$
\split
    * \ :&\ \omega=\omega|_{i_1,i_2,\cdots,i_p} d x^{i_1}
    d x^{i_2} \cdots  d x^{i_p} \\
     &\mapsto \sum_{j} \frac{\sqrt{g_M}}{(n-q)!}
     \omega|_{i_1,i_2,\cdots,i_p}
     \epsilon^{i_1,i_2,\cdots,i_p}_{\quad\quad\quad
      j_{p+1},j_{p+2},\cdots,j_n}
     d x^{j_{p+1}} d x^{j_{p+2}} \cdots d x^{j_n}
    ,\endsplit
$$
where $\epsilon_{i_1,i_2,\cdots,i_n}$ is a section of $\ZZ_{2 M}(U)$
$$
\epsilon_{i_1,i_2,\cdots,i_n} =\left\{ \matrix 
1 & \text{ if } (i_1,i_2,\cdots,i_n) 
\text{ is an even permutation of }(1,2,\cdots,n)\\
-1 & \text{ if }  (i_1,i_2,\cdots,i_n) 
\text{ is an odd permutation of }(1,2,\cdots,n)\\
0 & \text{ otherwise } \endmatrix \right. ,
$$
$$
 \epsilon^{i_1,i_2,\cdots,i_p}_{\quad\quad\quad j_{p+1},
 j_{p+2},\cdots,j_n}
 :=\tg_M^{i_1,j_1}\tg_M^{i_2,j_2}\cdots \tg_M^{i_p,j_p}
  \epsilon_{j_1,j_2,\cdots,j_p, j_{p+1},j_{p+2},\cdots,j_n}.
$$
Here we will express $d x^{i_1}  d x^{i_2} \cdots d x^{i_p}$ by
$d x^{i_1}\wedge  d x^{i_2}\wedge \cdots  \wedge d x^{i_p}$ 
for abbreviation.
}

\fvskip
If necessary, we can introduce an orthonormal system [Mal1] and
then the definition of the Hodge star operator becomes simpler.

\subheading{Proposition 3.3}

{\it 
$$
     \epsilon^{i_1,i_2,\cdots,i_n}=
    g_M^{-1}  \epsilon_{i_{1},i_{2},\cdots,j_n}.
$$
}

\demo{Proof} Direct computation gives the result. \qed \enddemo

\subheading{Proposition 3.4}

{\it 
\roster
\item
Let $w_M$ is the volume form of $(\Omega_M,\fg_M)$ and then
$$
    *1_M =w_M\in \Omega^n(M), \quad 
    *w_M = 1_M.
$$

\item  For $P \in \DM $, there exists $f \in \CM$ such that
$F_0(w_M\cdot P) \equiv f w_M$.

\endroster

}
\demo{Proof} From definition, (1) is obvious. Direct computations give (2).
\qed\enddemo

\subheading{Definition 3.5 (Adjoint Operator)}

{\it For  an open set $U$ of $M$ and $P \in \DM^{\CC}(U)$,
we define the adjoint operator $P^\dagger$ as
$$
    P^\dagger :=* (\overline{*1_M\cdot P})\equiv 
    * (\overline{w_M\cdot P}),
$$
where $ \overline{\ }$ means the complex conjugate.
If $P^\dagger= P$, we will call it self-adjoint operator.
}

\subheading{Remark 3.6}

We note that this adjoint operator can be expressed by the extension
operation [Bjo, TH].

\subheading{Definition 3.7}
{\it 

 The Beltrami-Laplace operator $\Delta$ of the Riemannian module $(M,\fg_M)$
 is defined by
$$
    \Delta_M = *d*d.
$$
}

\subheading{Proposition 3.8}
{\it 

There is a local system $\LL_M^n$  such that 
$\Theta_M(U) =\Gamma(U,\CM \otimes_{\RR_M}\LL_M^n)$ or
there is a basis $\{\theta_i\}_{i=1,\cdots,n}$ of $\LL_M^n(U)$,
$
\Theta_M(U)= \oplus_{i=1} \CM(U)\theta_i
$ for an open set $U$ in $M$.
}

\demo{Proof} 
From proposition 2.2, we can suppose that $\{x^i, \partial_i\}_{1\leq i\leq
n}$ is a local coordinate system for an open set $U$ of $M$ such that 
$[\partial_i, x^j]=\delta_i^j$. Then   $c_i \partial_i$ ($1\leq i\leq n, c_i
\in \RR_M(U) $) is an element of $\LL_M^n(U)$
and $\xi  = \sum_i a_i \partial_i \in \Theta_M(U)$ $a_i\in \CM(U)$.
\qed
\enddemo

Inversely,  for $\xi  = \sum_i a_i \partial_i \in \Theta_M(U)$, $a_i\in \CM(U)
$
and $\partial_i \in \LL_M^n$, we can find a local coordinate system $y^i$ 
for $V \subset (\cup_i \roman{supp}(a_i))\cap U$,
such that $\xi = \partial/ \partial y^j$,
where $\roman{supp}$ means the support of $a_i$.
This can be proved by solving the differential equation
$a_i = \partial y^j/\partial x^i$.

\subheading{Proposition 3.9}

{\it 
For a local coordinate $(x_1, \cdots, x_n)$ and inverse matrix $\tg_M^{ij}$ of
its metric  $g_{M ij}$ and $\tg_M^{ij}\equiv \tfg_M( d x^i, d x^j)$, the 
Beltrami-Laplace operator is expressed by $\Delta_M$ as $\Delta_M= g_M^{-1/2}
\partial_i g_M^{1/2} \tg_M^{i j } \partial_j$.
}

\demo{Proof}
For a local coordinate system, the volume form $w_M$ is expressed by
$w_M$ $= $ $g_M^{1/2}$ $ dx^1 \cdots d x^n$.
\qed\enddemo

Here we will define the Schr\"odinger system as a $\DMod$-module
related to the free  Schr\"odinger equation $-\Delta_M \psi =0$.

\subheading{Definition 3.10 (Schr\"odinger System)}
{\it 

We will define
the Schr\"odinger system $\SS_M$  by the exact sequence,
$$
    0 \longrightarrow \DM \overset
     -\Delta_M\to{\longrightarrow} \DM \longrightarrow 
    \SS_M \longrightarrow 0.
$$
}

As mentioned in the \S 1 introduction, we will give a lemma on
the Weyl algebra and thus introduce half form, wave system
and their related anti-self-adjoint differential operators as follows.

\subheading{Lemma 3.11}

{\it There exist elements $\{\partial_i\}_{1\leq i\leq n}$ $ \in  \LL_M^n$ 
such that they are of generators $\{ x^i,\partial_i\}_{1\leq i\leq n}$
 of Weyl algebra, $[x^i,\partial_j] =\delta^i_j$ 
 and  $(\sqrt{-1}\partial_i)^\dagger \neq \sqrt{-1}\partial_i$.}

\subheading{Definition 3.12 (Half Form)}

{\it We will define a sheaf of half form $\sqrt{\Omega^n_M}$ as
$\CM$-module whose section
over an open set $U$ of $M$ is given by $f \otimes_{ \RR_M} \sqrt{ w_M} $,
where $f\in \CM(U)$ and $\sqrt{ w_M}$ is defined as follows
\roster

\item $w_M := (\sqrt{ w_M})^2 = \sqrt{ w_M}\sqrt{ w_M}$,

\item $*\sqrt{ w_M}=\sqrt{ w_M}$,

\item For $\theta \in \Theta_M(U)$, the left handed action is given by
$$
    \theta(\sqrt{w_M}):=F_0(\theta\sqrt{w_M})= \frac{1}{2}f \sqrt{w_M},
$$
if  $F_0(\Lie_\theta(w_M) ) \equiv f w_M$.

\item The right handed action is
 $\sqrt{w_M}\cdot \theta = - \theta \sqrt{w_M}$.

\item For $\alpha \in \CM$,  
$ \theta \alpha \sqrt{w_M}= \theta(\alpha) \sqrt{w_M}+\alpha  \theta 
\sqrt{w_M}$. \endroster

We will denote $\CC_M \otimes_{\RR_M} \sqrt{\Omega_M^n}$
 by $\sqrt{\Omega_M^{\CC\  n  }}$.
}

If $M$ is a Riemannian surface with complex dimension one,
$\sqrt{\Omega_M^{\CC\  n  }}$ is essentially the same as
the prime form.
Further $\sqrt{w_M}$
essentially  appears in calculation of the gravitational or lorentzian
 anomaly in the elementary particle physics.
Accordingly it  is a natural variable.

\subheading{Definition 3.13 (Wave System $\sqrt{\CMC}$)}
{\it

A $\CMC$-module $\sqrt{\CMC}$, called wave system, is defined by
a closed presheaf
of the bilinear morphism; for $\psi_{a} \in \sqrt{\CMC}(U)$ ($a=1,2$), 
$*\psi_1 \cdot \psi_2\in \Omega_M^{\CC\  n  }(U)$.
}

\vskip 0.5 cm

We have the relation,
$$
    *\sqrt{\CMC}(U)\oplus \sqrt{\CMC}(U) \approx 
    \Omega_M^{\CC\  n  }(U) \approx  \CMC\cdot  w_M(U).
$$
We note that $\sqrt{\CMC}$ is isomorphic to $\CMC$ itself  because
basic field of $\CMC$ is a complete field $\CC$.

\subheading{Definition   3.14 (Anti-Self-Adjoint Connection 
$\tnablaM_\theta$ of Wave System)}
{\it

Let us define  an anti-self-adjoint connection $\tnablaM_\theta \in
 \SHom_{\CM}( \LL_M^n, \III(\sqrt{\CMC}) )$ by
 local relations for $ \theta\in \LL_M^n(U)$ and 
 $\sqrt{\omega_M}\in \sqrt{\Omega^n_M}(U)$,
$$
    \tnablaM_\theta:= *(\sqrt{w_M} \theta \sqrt{w_M}).
$$
Here we will denote a sheaf of set of $\tnablaM_\theta$ by
$\XiM$.
}

\subheading{Proposition  3.15 (Local Expression of $\tnablaM_\theta$)}
{\it

For a local coordinate system $\{x^i, \partial_i\}_{1\leq i\leq n}$
such that $[\partial_i, x^j]=\delta_i^j$,
 $\partial_i$ ($1\leq i\leq n$) is an element of $\LL_M^n$, and
 $w_M$ is expressed by 
$$
    w_M=\sqrt{g_M} dx_1 d x_2 \cdots dx_n,
$$
we have the local expression of $\tnablaM_{\theta}$ of $\theta=\partial_i$ cas
e,$$
\tnablaM_{\partial_i} = \partial_i + \frac{1}{4} \partial_i \log g_M= 
\root{4}\of{{g_M}^{-1}} \partial_i\root{4}\of{g_M} .
$$
}

\subheading{Corollary 3.16 (Local Expression of $\tnablaM_\theta$)}
{\it

For the local chart in Proposition 3.15 and $ \partial_{\alpha_1},
\partial_{\alpha_2},\cdots \partial_{\alpha_a}, \in \LL_M^n(U)$, 
$$ 
  \tnablaM_{\partial_{\alpha_1}}\tnablaM_{\partial_{\alpha_2}}\cdots 
  \tnablaM_{\partial_{\alpha_n}} =\root{4}\of{{g_M}^{-1}} 
  \partial_{\alpha_1}\partial_{\alpha_2}\cdots 
  \partial_{\alpha_a}\partial_i\root{4}\of{g_M} 
$$
}

\subheading{Proposition  3.17 (Anti-Self-Adjoint 
Connection $\tnablaM_\theta$ )}
{\it

\roster

For an open set $U$ of $M$, $\alpha \in \CM(U)$ and 
$\psi\in \sqrt{\CMC}(U)$,
following holds.

\item 
$$
    \tnablaM_{\alpha\theta}=\alpha\tnablaM_\theta.
$$

\item
$$
    \tnablaM_\theta(\alpha \psi)=\theta(\alpha) 
    \psi + \alpha \tnablaM_\theta( \psi).
$$

\item
$$
    [\tnablaM_\theta, \tnablaM_{\theta'}]=\tnablaM_{[\theta,\theta']}.
$$
\endroster
$\tnablaM_\theta$ is an integrable connection if it is a global section.
}

\demo{Proof}:
Hence, (1) and (2) are obvious. For $[\partial_i, \partial_j]=0$, 
$[\partial_i,\partial_j]g_M =0$. We obtain (4). \qed
\enddemo

\subheading{Definition  3.18 (Momentum Operator $\frak p_\theta$)}
{\it

Let $U$ is an open set of $M$.
\roster

\item
Let us define a momentum operator $\frak p_\theta $ of $\theta \in \LL_M^n(U)$
,
which consists of  an integrable connection $\tnablaM_{\theta} \in \XiM( 
\LL_M^n)(U)$ for $\theta \in \LL_M^n(U)$,
$$
    \frak p_\theta:= \sqrt{-1}\tnablaM_\theta.
$$

\item 
$ \Cal P_M$ is $\CC_M$-module generated by $\frak p_\theta$ for 
$\theta \in \LL_M^n(U)$, {\it i.e.}, $ \Cal P_M\equiv \sqrt{-1}\XiM$.

\endroster
}

\subheading{Proposition 3.19 (Momentum Operator $\frak p_\theta$)}
{\it
\roster

\item 
$\frak p_\theta$ of $\theta \in \LL_M^n$ is self-adjoint, {\it i.e.},
$\frak p_\theta^\dagger = \frak p_\theta$.

\item 
The Beltrami-Laplace operator $\Delta_M $ is locally expressed by
$$
\split 
\Delta_M =
\tnablaM_{i} \tg_M^{i j } \tnablaM_{j}&
 + \frac{1}{4} (\partial_j \log g_M)( \partial_i \tg_M^{i j})\\
&+ \frac{1}{4} \tg_M^{ij } ( \partial_i \partial_j\log g_M)+
\frac{1}{16} \tg_M^{ij }(\partial_i \log g_M)(\partial_j \log g_M).
\endsplit
$$
and $\tnablaM_{i} \tg_M^{i j } \tnablaM_{j}
=-\frak p_i\tg_M^{i j }\frak p_j$.
Here we have used notations $\tnablaM_i := \tnablaM_{\partial_i}$ and 
$\frak p_i:= \frak p_{\partial_i}$.

\endroster
}

\demo{Proof}:
$w_M\cdot \tnablaM_{\partial_1}= -\Lie_{\partial_1} w_M\ 
+\dfrac{1}{4}\partial_i \log g_M 
 w_M = w_M(- \theta - \dfrac{1}{4}\partial_i
\log g_M) $. From 3.15, (1) is obvious. By direct computation, we obtain 
$\tnablaM_{i} \tg_M^{i j } \tnablaM_{j}$.
\qed
\enddemo

\tvskip

\subheading{Corollary 3.20 }
{\it

$\Delta_M \equiv \frak p_i \tg_M^{i j } \frak p_i$
modulo $F_0( \DM)$.

}

\subheading{Definition 3.21 (Momentum Operator $\frak q_\theta$)}
{\it
\roster

Let $U$ is an open set of $M$.

\item
Let us define a momentum operator $\frak q_\theta $ of $\theta \in \LL_M^n(U)$
,
which consists of  an integrable connection $\theta$ of
$\Gamma\left(U, \SHom_{\RR_M}(\LL_M^n, \III(\sqrt{\Omega_M^{\CC\ n}}) )\right)
$,$$
    \frak q_\theta:= \sqrt{-1}\partial_\theta,
$$
where for $\alpha \in \CM(U)$ and $f \sqrt{w_M}\in \sqrt{\Omega_M^{\CC\ n}}(U)
$,$$
\split
    \theta(\alpha f \sqrt{w_M})&=\theta(\alpha) 
    f\sqrt{w_M}  + \alpha \theta( f \sqrt{w_M} ),\\
    &=\theta(\alpha f) \sqrt{w_M}  + \alpha f\theta( \sqrt{w_M} ),\\
\endsplit
$$

\item
$\Cal Q_M$ is a $\CC_M$-module sheaf whose local section 
is $\frak q_\theta$.

\endroster
}

\subheading{Remark 3.22 (Momentum Operator $\frak p_\theta$)}
{\it

Noting $*$ for $\sqrt{\Omega_M}$ is identity operator
 and we have the relation,
$$
    *\sqrt{\Omega_M^{\CC\  n  }} \oplus \sqrt{\Omega_M^{\CC\  n  }} 
    \approx\Omega_M^{\CC\  n  } \approx  \CMC\cdot  w_M.
$$
We should modify the definition of $\dagger$ and 
let $\frak q_\theta$ be  self-adjoint;
$
      \frak q_\theta^\dagger = \frak q_\theta.
$
}

\subheading{Proposition 3.23 }

{\it
\roster 

Following  categories are isomorphic.

\item 
Let $\Cal A$ be a category  whose objects are $\CM$-modules $\sqrt{\CMC}$ and 
morphisms are elements of  a differential ring sheaf $\DMod_{\Cal P_M}^\CC$
generated by $\Cal P_M$
with coefficient  $\sqrt{\CMC}\approx{\CMC}$.

\item Let $\Cal B$ be a category  whose objects are $\CM$-modules 
$\sqrt{\Omega_M^{\CC\ n}}$ and morphisms are elements of a differential ring
sheaf $\DMod_{\Cal Q_M}^\CC$ generated by $\Cal Q_M$ with coefficients ${\CMC}
$.
\endroster
In other words, there is an equivalent functor,
$\xi_M:(\sqrt{\CMC},\DMod_{\Cal P_M}^\CC) \to (\sqrt{\Omega_M^{\CC 
n}},\DMod_{\Cal Q_M}^\CC)$.
}

\demo{Proof}
Let $\xi_M:(\sqrt{\CMC},\tnablaM_\theta) \mapsto (\sqrt{\Omega_M^{\CC 
n}},\theta)$ for $\theta$ in $\LL_M^n(U)$. From the assumption, there is no 
point $p$ in $M$ such that $g_M(p)=0$. We denote $*1_p = \sqrt{g_M} dx^1\cdots
dx^n=:\sqrt{g_M} d^n x$. Noting the proposition 3.15, for sections of $\psi \i
n
\sqrt{\CMC}(U)$ and $\theta \in \LL_M^n$, the correspondences $\xi_M (\psi) =
\sqrt{d^n x}\root{4}\of{g_M}\psi $ and $\xi_M(\tnablaM_\theta)=$ $\xi_M (  
\root{4}\of{{{g_M}^{-1}}}\theta \root{4}\of{g_M})= \theta$ are bijection as a
set due to non-degeneracy of $g_M$. Thus we will check the commutativity of 
$\xi_M $ and actions of derivative, {\it i.e.}, $\xi_M \tnablaM_i = \partial_i
\xi_M$ for $\theta\equiv \partial_i$ case; $\xi_M (\tnablaM_i \psi) $ $= 
\sqrt{d^n x} \root{4}\of{g}\tnablaM_i \psi $ $ =\sqrt{d^n x} \root{4}\of{g_M}
(\root{4}\of{{g_M}^{-1}}\partial_i  \root{4}\of{g_M}) \psi  $ $= \sqrt{d^n 
x}\partial_i  \root{4}\of{g_M}\psi = \partial_i \xi_M(\psi)$. We completely 
prove it.
\qed\enddemo

\subheading{Lemma 3.24 }
{\it

\roster
Let $U$ be an open set of $M$.

\item There is an equivalent functor $\tilde \sigma_M$ from $\DMod_M^\CC$ to
$\DMod_{\Cal Q_M}^\CC$. For $P \in \DMod_M^\CC$, $g_M^{1/4} P 
g_M^{-1/4}\in\DMod_{\Cal Q_M}^\CC$.

\item There is an equivalent functor $\overline{\sigma}_M$ from $\DMod_M^\CC$
to $\Cal D_{\Cal P_M}^\CC$; $\overline{\sigma}_M:=\xi_M^{-1}\circ \tilde 
\sigma_M$. For an element of $P \in \DMod_M^\CC(U)$, $\overline{ \sigma}(P)=P$
.

\item
$P'\in \DMod_{\Cal Q_M}^\CC(U)$ locally has a standard form representation,
$$
           P' = \sum_{\alpha \in \Bbb N^n} a_\alpha\partial^\alpha .
$$

\item
$P'\in \DMod_{\Cal P_M}^\CC(U)$ locally has a standard form representation,
$$
           P' = \sum_{\alpha \in \Bbb N^n} a_\alpha\tnablaM^\alpha .
$$

\item The functors $\tilde \sigma_M$ and $\overline{\sigma}_M$ are extended to
functors of the category of left $\DMod_M^\CC$-modules and that of 
left $\Cal D_{\Cal Q_M}^\CC$-modules
or left $\Cal D_{\Cal P_M}^\CC$-modules.

\endroster
}

\demo{Proof}
From the propositions 3.15 and 3.23 give the correspondences (1) and (2), whic
h
are essentially the same as the hermitianization known in the Sturm-Liouville
operator [Arf] as showed in the introduction. For an element $\partial_i$ of
$\Theta_M(U)$, we have a local expression, 
$\overline{\sigma}_M(\partial_i)=\tnablaM_i - (\partial_i \log g_M)/4 =$ 
$g_M^{1/4} [g_M^{-1/4} \partial_i g_M^{1/4}] g_M^{-1/4}=\partial_i$. The 
expression of (3) can be obtained by expanding $g_M^{-1/4}\partial^\alpha 
g_M^{-1/4}=\partial^\alpha + \cdots$. (4) is guaranteed from (3) and 
$\xi^{-1}$. Next we will show (5). For a quotient space $\Cal M= 
\DMod_M^\CC/\DMod_M^\CC R$ of $R \in \DMod_M^\CC$ and $P_1$, $P_2$, $Q \in 
\DMod_M^\CC$ such that $P_1 - P_2=Q R$, $\tilde \sigma_M(P_1-P_2) $ $=g_M^{1/4
}
Q g_M^{-1/4}g_M^{1/4}R g_M^{-1/4}$ $=\tilde \sigma_M(Q) \tilde \sigma_M(R)$.
Thus $\tilde \sigma_M(\Cal M)$ can be defined as $\tilde \sigma_M(M):= 
\DMod_{\Cal Q_M}^\CC/\DMod_{\Cal Q_M}^\CC \tilde \sigma_M(R)$. Similarly we ca
n
naturally define $ \tilde \sigma_M$ for general coherent module.
\qed \enddemo

\tvskip 

Due to the properties of $\overline{\sigma}_M$, we can mix 
$\overline{\sigma}_M(\DMod_M^\CC)$ and $\DMod_M^\CC$. In fact we did not 
discriminate them  in the introduction and will not in [II]. However in this
article, in order to see the action of $\overline{\sigma}_M$, we will 
explicitly express it.

\subheading{Corollary 3.25 }
{\it

\roster

\item
In the category $\Cal B$, we have local expression
of $\tilde \sigma_M(\Delta_M)=g_M^{1/4} \Delta_M g_M^{-1/4}$,
$$
\split
\tilde\sigma_M(\Delta_M)&=\partial_i \tg_M^{i j } \partial_i 
+ \frac{1}{4} (\partial_j \log g_M)
( \partial_i \tg_M^{i j})\\
&+ \frac{1}{4} \tg_M^{ij } ( \partial_i \partial_j\log g_M)+
\frac{1}{16} \tg_M^{ij }(\partial_i \log g_M)(\partial_j \log g_M).
\endsplit
$$

\item
$\overline{\sigma}_M(\Delta_M)$ agrees with proposition 3.19 (2),
$\Delta_M=\overline{\sigma}_M(\Delta_M)$.

\endroster
}

\fvskip

{\centerline{\bf{\S 4. Schr\"odinger  Operators 
in a submanifold $S \hookrightarrow\EE^n$}}
\fvskip

Let $S$ be a $k$-dimensional real  analytic compact submanifold of $\EE^n$, 
$\iota_S: S \hookrightarrow \EE^n$ and $T_S$ be its associated  tubular 
neighborhood; $i_S :S \hookrightarrow T_S$  and $ i_{T_S}: T_S \hookrightarrow
\EE^n$ such that $\iota_S \equiv i_{T_S} \circ i_S$. $T_S$ has a projection 
$\pi_{T_S} : T_S \to S$. Since $\EE^n$ has a natural metric $\fg_{\EE^n}$,
$T_S$ and $S$ have Riemannian modules $(i_{T_S}^{-1} \Theta_{\EE^n}, 
i_{T_S}^{-1} \fg_{\EE^n})$ and $(\Theta_S, \iota_{S}^{-1} \fg_{\EE^n})$. 
Further we put $i_S^*(\XiT(\Theta_{T_S}^{\perp}))
:=\CS \otimes_{i_S^{-1} \CTS}i_S^{-1}\XiT(\Theta_{T_S}^{\perp})$.

The Schr\"odinger system ${\SS}_{T_S}$ of $T_S$ is given as 
${\SS}_{T_S\to \EE^n}:=i_{T_S}^{*}{\SS}_{\EE^n}$ 
$:=\DMod_{T_S\to \EE^n}\otimes_{i_{T_S}^{-1}
{\DMod}_{\EE^n}}i_{T_S}^{-1}{\SS}_{\EE^n}$
and  
$$
  \overline{\SS}_{S \to T_S} 
  :=i_{S}^{*} \overline{\sigma}_{T_S}({\SS}_{T_S})
  :=\overline{\sigma}_S(\DMod_{S})\otimes_{i_{S}^{-1}
  \overline{\sigma}_{T_S}({\DMod}_{T_S})}i_{S}^{-1} 
   \overline{\sigma}_{T_S}({\SS}_{T_S}).
$$
Further we will use the notations in \S 2 and \S 3 and 
 will not neglect the action of $\overline{\sigma}$'s.

\subheading{Proposition 4.1}
{\it

\roster 
\item
For $i_{S}^{*}{\SS}_{T_S}$ 
$:=\DMod_{S}\otimes_{i_{S}^{-1}{\DMod}_{T_S}}i_{S}^{-1}{\SS}_{T_S}$,
$$
i_{S}^{*}{\SS}_{T_S}= \iota_{S}^{*}{\SS}_{\EE^n}
:= \DMod_{S}\otimes_{\iota_{S}^{-1}{\DMod}_{\EE^n}}\iota_{S}^{-1}{\SS}_{\EE^n}
.
$$

\item
 Let the anti-self-adjoint connection 
 $\tnablaS_\dalpha\in\Gamma(U,i_S^*\XiT(\Theta_{T_S}^{\perp}))$ 
 for  an open set $U$ in $S$.
There is an injective endomorphism of 
the Schr\"odinger system $\overline{\SS}_{S \to T_S} $,
$$
\eta^\conf_\dalpha 
:\overline{\SS}_{S \to T_S}  \to \overline{\SS}_{S \to T_S} ,
$$
for $P \in \Gamma(U,\overline{\SS}_{S \to T_S} ) $,
$
       \eta^\conf_\dalpha(P) =  P \tnablaS_\dalpha \in
        \Gamma(U,\overline{\SS}_{S \to T_S} ).
$
Then we have a submodule of $\overline{\SS}_{S \to T_S} $, 
$$
\eta^\conf:(\overline{\SS}_{S \to T_S} )^{n-k}
\to \sum_{\dalpha=k+1}^n\overline{\SS}_{S \to T_S}  
\tnablaS_\dalpha \subset \overline{\SS}_{S \to T_S} .
$$
}

\demo{Proof}
Due to the relation  $\iota_S \equiv i_{T_S} \circ i_S$, (1) is obvious. From
lemma 3.24, we choose elements $P_1$, $P_2$ and $Q$ in $(i_{S}^{*} 
\overline{\sigma}_{T_S}({\DMod}_{T_S}^\CC))$ such that $P_1\equiv P_2\in 
\overline{\SS}_{S \to T_S} $ {\it i.e.}, $P_1-P_2 = Q 
[i_{S}^{*}\overline{\sigma}_{T_S}(\Delta_{T_S})]$. 
Noting the relation
$$
	[i_{S}^{*}\overline{\sigma}_{T_S}(\Delta_{T_S}), \tnablaS_\dalpha]=0,
$$
$Q [i_{S}^{*}\overline{\sigma}_{T_S}(\Delta_{T_S})]\tnablaS_\dalpha$
$=Q \tnablaS_\dalpha[i_{S}^{*}\overline{\sigma}_{T_S}(\Delta_{T_S})] \equiv 0 
$
in $ \overline{\SS}_{S \to T_S} $. Thus it is injective.
\qed \enddemo

\subheading{Definition 4.2}
{\it
\roster

\item
We will define a coherent $i_{S}^{*} 
\overline{\sigma}_{T_S}(\DMod_{T_S})$-module ${\bar \SS}_{S}^{\parallel}$ in
$T_S$ by the exact sequence,
$$
         (\overline{\SS}_{S \to T_S} )^{n-k}
     \overset{\eta^\conf}\to{\longrightarrow}
     \overline{\SS}_{S \to T_S} \longrightarrow
   {\overline \SS}_{S\to T_S}^{\parallel}\longrightarrow 0.
$$

\item
We will define a coherent $\DMod_{S}$-module  by
$$
    {\SS}_{S\hookrightarrow \EE^n}^{}:=\overline{\sigma}_{S}^{-1} 
    {\overline \SS}_{S\to T_S}^{\parallel}.
$$
Let us call it  submanifold Schr\"odinger system.

\item When the submanifold Schr\"odinger system ${\SS}_{S\hookrightarrow 
\EE^n}$ is decomposed by the exact sequence,
$$
0 \longrightarrow \DMod_S \overset {{-\Delta}_{S\hookrightarrow \EE^n}}
\to{\longrightarrow} \DMod_S \longrightarrow 
    {\SS}_{S\hookrightarrow \EE^n}\longrightarrow 0,
$$
where ${\Delta}_{S\hookrightarrow \EE^n}-{\Delta}_{S} \in F_0 \DMod_S$,  
we will call
${\Delta}_{S\hookrightarrow \EE^n}$ 
 the submanifold Schr\"odinger operator.
\endroster
}

\subheading{Proposition 4.3}
{\it

\roster

\item ${\overline \SS}_{S\hookrightarrow \EE^n}$ is uniquely determined.

\item
The definitions in 4.2  are naturally  extended to a submanifold immersed
in $\EE^n$.
\endroster
}

\demo{Proof}
From the definition which does not depend upon the coordinate system,
 (1) is obvious. 
When we construct ${\Delta}_{S\hookrightarrow \EE^n}$ and others,
we used only local date. Hence (2) is also obvious. \qed \enddemo

Now we will state our main theorem in this article.

\subheading{Theorem 4.4}
\roster

\item $k=1$ and $n=3$ case, $S$ is curve $C$,

$$
   -\Delta_{C \hookrightarrow \EE^3} 
   = -\partial_s^2 - \frac{1}{4} |\kappa_\CC|^2,
$$
where $s$ is the arclength of the curve $C$, $\kappa_\CC$ is the complex 
curvature of $C$, defined by $\kappa_\CC = \kappa(s) \exp\left(\sqrt{-1} 
\dsize\int^s \tau ds\right)$ using the Frenet-Serret curvature $\kappa$ and
torsion $\tau$.

\item $k=2$ and $n=3$ case, $S$ is a conformal surface,

$$
          -\Delta_{S \hookrightarrow \EE^2} = -\Delta_S - (H^2 -K),
$$                                         
where 
$H$ and $K$ are the mean and Gauss curvatures.

\endroster

These operators agree with the operators obtained by Jensen and Koppe [JK]
and da Costa [dC].
Investigation of these equations might mean the properties of these 
submanifolds. Similar attempt was done [HL] for $-\Delta_S -2 H^2$ but these
operators are more natural because they are related to Frenet-Serret and 
generalized Weierstrass relations [II].

\fvskip

In order to prove theorem 4.4, we will set up the language to express the 
submanifold system. We will note that these concepts of differential 
geometry, such as the Christoffel symbol, curvature and so on, are translated
to language in sheaf theory by Mallios [Mal1,2]. Thus although we will use 
them in classical ways, they could be written more abstractly if one prefers.

An affine vector in $\EE^n$ is given by $(Y^1,Y^2, \cdots, Y^n)$ 
as the Cartesian coordinate system and 
in its tangent space  $T_p\EE^n$, the  bases are
$\partial_i:=\partial/\partial Y^i$, $i=1,2,\cdots,n$,
 $<\partial_i,\dd Y^j>=\delta_i^j$. 
 Here Latin indices $i,j,k,$ are for the Cartesian coordinate of $\EE^n$. 
 $\Theta_{\EE^n} = \CEE \otimes_{\RR_{\EE^n}} \LL_{\EE^n}^n$.
The Riemannian  metric $\fg_{\EE^n}$ in the euclidean space is 
given as
$$
    \fg_{\EE^n}=\delta_{i,j} \dd Y^i \otimes \dd Y^j.
$$

Let the equations $Q^a(Y^1,Y^2,\cdots,Y^n)=0$, ($a=1,\cdots,d\equiv n-k$ expre
ssa surface of $S$; the Pfaffian is expressed by $\dd Q^a=0$'s and  by Frobeni
us 
integrable theorem, there are vector fields given by the bases 
$\partial_\alpha:=\partial/\partial s^\alpha$ which are satisfied with
 $$
     <\partial_\alpha,\dd Q^\dalpha>=0 ,\quad \text{for } 
     \forall \alpha \  \forall \dalpha.
 $$ 
Hence  the local coordinate of the 
submanifold is given as $(s^1,s^2,\cdots,s^k)$ or $(s^\alpha)$.
Let us employ the conventions
 that the beginning of the Greek ($ \alpha$, $\beta$, $\gamma$, $\cdots$)
 runs from 1 to $k$ and it with dot 
 ($ \dalpha$, $\dbeta$, $\dgamma$, $\cdots$) runs form $k+1$ to $n$.

\subheading{Notation 4.5}

{\it 

\roster

\item A point $p$ in $T_S$  is expressed by the  local coordinate 
$(u^\mu):=(s^1,s^2,\cdots,s^k,q^{k+1},\cdots,$ $q^n)$, $\mu=1,2,\cdots,n$ wher
e
$(s^1, \cdots,s^k)$ is a local coordinate of $\pi_{T_S} p$; We assume that
the beginning of the Greek 
($ \alpha$, $\beta$, $\gamma$, $\cdots$) runs from
1 to $k$ and they with dot ($ \dalpha$, $\dbeta$, $\dgamma$, $\cdots$) runs
form $k+1$ to $n$).

\item Let
$(u^\mu)=(s^\alpha,q^\dalpha)$, where the middle part of the Greek ($ \mu$,
$\nu$, $\lambda$, $\cdots$) run from $1$ to $n$.

\item $\bee_\alpha:=\partial_{\alpha}:=\partial/ \partial s^{\alpha}$ is a
base of $ \Theta_S(U)$. For $\bee_\alpha \in \iota_S^{-1} \Theta_{\EE^n}(U)$, 
$\bee_\alpha$ is expressed by $ \bee_\alpha=e^i_{\ \alpha}\partial_i$.

\item Using $<\bee_\alpha, \bee^\beta>=\delta^\beta_\alpha$, $\bee^\beta 
\equiv \dd s^\beta\in \Omega_{S}^1(U)$.

\item $\bee_\dalpha:=\partial_{\dalpha}:=\partial/ \partial q^{\dalpha}$ is
a base of $ \Theta_{S^\perp}(U)$. For $\bee_\dalpha \in \iota_{S}^{-1} 
\Theta_{\EE^n}(U)$, $\bee_\dalpha$ is expressed by 
$ \bee_\dalpha=e^i_{\ \dalpha}\partial_i$.

\item Let $
    (\iota_S^{-1} \fg_{\EE^n})(\bee_\mu,\bee_\nu)=\delta_{i,j} 
    e_{\ \mu}^{ i} e_{\ \nu}^{j} =g_{\EE^n \mu,\nu}$.
Then we have induced metric   $\fg_{ T_S}:=i_{T_S}\fg_{\EE^n} $ and 
$\fg_{S}:= \iota_S^{-1} \fg_{\EE^n}$. Then we will express $\fg_{S}$ by
the relation, 
$$   
    \fg_{S}= \iota_S^{-1}( \delta_{i,j}\dd (x^i)\otimes \dd (x^j)) 
     =g_{S,\alpha\beta} \dd s^\alpha \otimes \dd s^\beta , 
$$
where
$$
    g_{S,\alpha\beta}:=\fg_{S}(\partial_\alpha,
    \partial_\beta):= \iota_S^{-1} \fg_{\EE^n}
     ( \bee_\alpha,\bee_\beta).
$$

\endroster
}

}

\subheading{Proposition 4.6}

{\it 
Let $\overline{\DMod}_{S \to T_S}$  
$:=\CS \otimes _{i_S^{-1}\CT} (i_S^{-1}\overline{\sigma}_{T_S}(\DMod_{T_S})$.
\roster

\item
For an element $P$ of left $\overline{\DMod}_{S \to T_S}$-module 
$\overline{\SS}_{S\to T_S}$, there exists $Q\in 
\overline{\sigma}_{T_S}{\SS}_{T_S}$ such that 
$P = Q|_{q^\dalpha=0,\ {\dalpha= k+1,\cdots,n}}$ 
for a normal coordinate $(q^\dalpha,\ \dalpha = k+1,\cdots,n)$.

\item
For an element $P$ of left $\DMod_{S}$-module 
$\overline{\SS}_{S\to T_S}^\parallel$,
there exists $Q\in\overline{\SS}_{S\to T_S}$ such that  
$P = Q|_{\tnablaT_\dalpha=0,\ {\dalpha = k+1,\cdots,n}}$.

\item ${\Delta}_{S\hookrightarrow \EE^n}$ is uniquely determinded.
\endroster
}

\demo{Proof} 
 (1) and (2) are obvious from the definition 4.2. Since we tuned 
 ${\Delta}_{S\hookrightarrow \EE^n}$ using ${\Delta}_{S}$, there is no 
 multiplicative freedom. From propositoon 4.3 (1), (3) is obvious.
  \qed \enddemo

\subheading{Proposition 4.7}

{\it An affine vector (coordinate) $\bold Y\equiv(Y^i)$ in $T_S \subset \EE^n$
is expressed by,
$$
     \bold Y= \bold X + \bee_\dalpha q^\dalpha,
$$
for a certain affine vector $\bold X$ of $S$. 
}

\demo{Proof} By setting $\bold X = \pi_{T_S} \bold Y$, it is obvious.\qed
\enddemo

For a case of immersion, this expression is not unique but locally unique.

\subheading{Proposition 4.8}

{\it 
For $U \subset T_S$, the induced metric of $T_S$ from $\EE^n$ has a direct
sum form,
$$
\fg_{T_S}:=i_{T_S}^{-1}\fg_{\EE^n}
=\fg_{T_S^\parallel}\oplus\fg_{T_S^\perp},
$$
where $\fg_{T_S^\perp}$ is trivial structure.
In local coordinate, 
$$
\fg_{T_S^\perp}= \delta_{\dalpha,\dbeta} 
\dd q^\dalpha\otimes \dd q^\dbeta,
\quad \fg_{T_S^\parallel}
     =g_{T_S\alpha\beta} \dd s^\alpha \otimes \dd s^\beta ,
$$
or for $g_{T_S\mu,\nu}:=\fg_{T_S}(\partial_\mu,\partial_\nu)$
$$
 g_{T_S\dalpha\dbeta}=\delta_{\dalpha\dbeta},
     \quad
     g_{T_S\dalpha\beta}=g_{T_S\alpha\dbeta}=0,
$$
where $\partial_\mu:=\partial/\partial u^\mu$.
}

\fvskip

In order to prove this proposition and to
give a concrete expression of $\fg_{S^\parallel}$
 using $\fg_{S}$ and $q^a$,
we will consider the intrinsic and the extrinsic properties of
 $S \subset \EE^n$
{\it e.g.}, the Weingarten map [E,G].

\fvskip

\subheading{ Proposition 4.9 (intrinsic properties)}

We can define the Riemannian connection consisting with this metric
$\fg_{ S}$ for $\theta,\xi \in \Theta_S$,
$$
    D_\theta \xi:= \theta \xi - \fg_{S}(\theta \xi, \bee_\dbeta) \bee_\dalpha
    \tg_{T_S}^{\dalpha,\dbeta}.
$$
where $\tg_{T_S}^{\dalpha,\dbeta}$ is the inverse matrix of 
$g_{T_S\dalpha,\dbeta}$.

\demo{Proof} See chapter 12 and 13 in [G]. \qed \enddemo

\subheading{ Proposition 4.10 (Weingarten map)}

{\it 
For a  base $\bee_\alpha$ of $\Theta_S$ and  $\bee_\dbeta \in \iota_S^{-1}
\Theta_{\EE^n}$, $\partial_\alpha\bee_\dbeta$ 
is an element of  $\iota_S^{-1}
\Theta_{\EE^n}$ and is expressed by 
$$
    \partial_\alpha \bee_\dbeta = \gamma_\dbeta(\partial_\alpha) = 
    \gamma^\dalpha_{\ \dbeta\alpha} \bee_{\dalpha}+\gamma^\beta_{\ 
    \dbeta\alpha} \bee_{\beta}
    ,
\quad
    \gamma_{\ \dbeta\alpha}^\mu :=\fg_{ S}( 
    \gamma_\dbeta(\partial_\alpha), \bee^\mu). 
$$
$-\gamma_\dbeta: \Theta_{S} \to \iota_S^{-1} \Theta_{\EE^n}$
 is called as the
Weingarten map. }
  
\demo{Proof} See chapter 12 and 13 in [G] and p.162-164 in [E]. \qed\enddemo

\subheading{ Proposition 4.11 (Second fundamental Form)}

{\it 

The second fundamental 
$
    \gamma^\dalpha_{\ \beta\alpha}:=\fg_{S}
    (\partial_\alpha \bee_\beta, \bee_\dbeta) 
$
is connected with the Weingarten map,
$$
    \gamma^\dalpha_{\ \beta\alpha}
    =-\fg_{S}(\bee_\beta,\gamma^\gamma_{\ \dalpha \alpha}
    \bee\gamma) , \quad
    \gamma^\dalpha_{\ \beta\alpha}
    =-g_{S, \beta\gamma} \gamma^{\gamma}_{\ \dalpha \alpha} .
$$
}

\demo{Proof} [E]
Due to $\fg_{S}(\bee_\alpha,\bee_\dbeta)=0$ and
 $\partial_\beta \fg_{S}(\bee_\alpha,\bee_\dbeta)=0$,
 we prove it. 
 \qed \enddemo

\fvskip

\subheading{Lemma 4.12}
{\it

There exist the normal vectors $\bee_\dalpha\in \Theta_{S^\perp}$
satisfied with,
$$
    \partial_\alpha {\bee}_{\dalpha}
    =\gamma_{\ \dalpha \alpha}^{\beta} 
    {\bee}_{\beta}.
$$
}

\demo {Proof} Let ($\gamma,\bee)$ in proposition 4.11 be rewrite  
$(\tilde \gamma,\tilde{\bee}$).
From the proposition 4.11, 
the derivative of a general normal orthonormal base 
$\tilde{\bee}_\dalpha$ is given as
$\partial_\alpha \tilde{\bee}_{\dalpha}
    = \tilde\gamma_{\ \dalpha \alpha}^{\beta} 
    \tilde{\bee}_{\beta}.
         + \tilde\gamma_{\ \dalpha \alpha}^{\dbeta}
         \tilde{\bee}_{\dbeta} $. 

From $\fg_{\EE^n}(\bee_\dalpha,\bee_\dbeta)=\delta_{a,b}$, for $\theta \in
\Theta_S(U)$, {\it i.e.}, $\theta=f^\alpha \partial_\alpha$
at $U\subset S$,
$
    \fg_{\EE^n}(\partial_\alpha\bee_\dbeta,\bee_\dalpha)
    =- \fg_{\EE^n}(\bee_\dbeta,\partial_\alpha\bee_\dalpha) 
$,
we have
$
    \tilde\gamma_{\ \dalpha \alpha}^{\dbeta}
    =- \tilde\gamma_{\ \dbeta \alpha}^{\dalpha},$ and $
      \tilde\gamma_{\ \dalpha \alpha}^{\dalpha}\equiv 0$
      (not summed over $\dalpha$).
In other words, there are $k (n-k)(n-k-1)/2$ degrees of freedom; 
$\tilde\gamma_{\ \dbeta\alpha}^{\dalpha}$ for $\alpha=1,\cdots,k$.
Thus we will employ an element $G$ of SO$(n-k)$ transformation so that 
$$
    (\partial_\alpha + \tilde\gamma_{\ \dbeta\alpha}^{\dalpha})=
        G^{-1}(\partial_\alpha)G.
$$
It is obvious that the solution of this differential 
equation locally exists, 
{\it e.g.},
$$
    \pmatrix \bee_\dalpha \\ \bee_\dbeta 
    \endpmatrix= G_{\dalpha\dbeta} 
      \pmatrix \tilde{\bee}_{\dalpha}\\ 
      \tilde{\bee}_{\dbeta} \endpmatrix,
$$
$$
G_{\dalpha\dbeta} 
      = \pmatrix \cos\theta &-\sin\theta \\
       \sin\theta &\cos\theta \endpmatrix,
\quad
    \theta := \int^{s^1} \dd s^1\ 
    \tilde\gamma_{\ \dbeta1}^{\dalpha}+\int^{s^2} \dd s^2\ 
    \tilde\gamma_{\ \dbeta2}^{\dalpha}+\cdots
    +\int^{s^k} \dd s^k\ 
    \tilde\gamma_{\ \dbeta k}^{\dalpha}.
$$
The topological structure of $\EE^n$ is simple and the normal
bundle of $S$ exists if the submanifold $S$ is not wild.
From the assumptions on $S$, there is no singularity in $S$.
Thus these solutions globally  exist.
This transformation is sometimes called as Hashimoto
 transformation.
\qed \enddemo

\fvskip

\subheading{ Lemma 4.13}
{\it

The moving frame of $T_S$,
$\bEE_\mu=\pi_{T_S}^{-1}(\partial_\mu)\in \Theta_{T_S}$
($\mu=1,\cdots,n$), 
is expressed by $\bEE_\mu=E_{\ \mu}^i \partial_i \in i_{T_S}^{-1} 
\Theta_{\EE^n}$ and
$$
    E^i_{\ \alpha} = e^i_{\ \alpha} + 
            q^\dalpha \gamma^\beta_{\ 
            \dalpha\alpha} e^i_{\ \beta}, 
            \quad  E^i_{\ \dalpha} = e^i_{\ \dalpha}.    
$$
}

\demo{Proof} Using proposition 4.7, direct computation leads this result.
\qed \enddemo

\subheading{Proof of Proposition 4.8}

 Lemma 4.12 and $\fg_{TS}=E_{\ \mu}^i E_{\ \mu}^i d x^i \otimes d x^i$
lead the result in proposition 4.8.
\qed 

\fvskip

Its inverse matrix is denoted by $ (E^\mu_{\ I})$.

\subheading{Corollary 4.14}

{\it 
\roster
\item The metric in $T_S$ is expressed as
$$
    \fg_{T_S^\parallel}=\fg_{S}+
    \fg^{(1)}_{TS} q^\dalpha+\fg^{(2)}_{TS} (q^\dalpha)^2  ,
$$
$$
g_{T_S^\parallel}(\partial_\alpha,\partial_\beta)
                        =
g_{S \alpha\beta}+
    [\gamma_{\ \dalpha\alpha}^\gamma g_{S\gamma\beta}+
    g_{S\alpha\gamma}\gamma_{\ \dalpha\beta}^\gamma]q^\dalpha
    +[\gamma_{\ \dalpha\alpha}^\delta g_{S\delta\gamma}
     \gamma_{\ \dbeta\beta}^\gamma]q^\dalpha q^\dbeta .
$$

\item $g_{T_S}:=\det_{n\times n}(g_{T_S,\mu.\nu})$ is $
g_{T_S}=\det_{k\times k}(g_{T_S^\parallel,\alpha.\beta})$ and 
$$
\split
    g_{T_S}=g_{S} &\Bigl\{1+2 \tr_{k\times k}
    (\gamma^{\alpha}_{\ \dalpha\beta})q^\dalpha\\
      &+\left[ 2 \tr_{k\times k}(\gamma^{\alpha}_{\ \dalpha\beta})
          \tr_{k\times k}(\gamma^{\alpha}_{\ \dbeta\beta})-
          \tr_{k\times k}(\gamma^{\delta}_{\ \dalpha\beta}
         \gamma^{\alpha}_{\ \dbeta\delta}) 
         \right] q^\dalpha q^\dbeta+
         \Cal O(q^\dalpha q^\dbeta q^\dgamma) + \cdots \Bigr\}.
\endsplit
$$
\endroster
}
\vskip 1.0 cm

\subheading{ Example 4.15}
\roster

\item
In coordinate,  for the case of $n=3, k=1$;
$$
    g_{T_S^\parallel}=g_{S} (1-|\kappa_\CC| q^3+|\kappa_\CC|^2 (q^3)^2)^2 ,
$$
where  $\kappa_\CC:=\gamma_{\ 12}^3 +\sqrt{-1}\gamma_{\ 13}^2 $ 
is the complex curvature of $C$. Here 
$\kappa_\CC $  $= \kappa(s) \exp\left(\sqrt{-1} \dsize\int^s
\tau ds\right)$ is also given by the Frenet-Serret 
curvature $\kappa$ and torsion $\tau$.

\item
In coordinate,  for the case of $n=3, k=2$;
$$
    g_{T_S^\parallel}=g_{S} (1-2H q^3+K (q^3)^2)^2,
$$
where $H:= \tr(-\gamma_{3\beta}^\alpha)/2$ is the mean curvature and 
$K:=\det(-\gamma_{3\beta}^\alpha)$ is the Gauss curvature.

\item 
In coordinate,  for the case of $n>3$ and $ k=2$;
$$ 
    g_{T_S}=g_S 
        (1+\roman{tr}_2(\gamma^\alpha_{\ 3\beta}) q^3
        +\roman{tr}_2(\gamma^\alpha_{\ 4\beta})q^4
                +K(q^3,q^4))^2.
$$
We will denote
$$
H_{\dalpha-2}:=-\frac{1}{2}\roman{tr}_2(\gamma^\alpha_{\ \dalpha\beta}).
$$
We can introduce the "complex mean curvature" for $n=4$ case,
$$
    H_c=H_1+ \ii H_2. 
$$

\endroster

\demo{Proof}
Since $g_{T_S^\parallel}=(\det(\partial_\alpha x^i))^2$, 
we obtain (1) after calculation of 
$\det(\partial_\alpha x^i)$ using the fact that for a $2 \times 2$
 matrix $A$, $\det(1+A)=1 + \tr A + \det A$. Similarly we have (2) and (3).
\qed \enddemo
\fvskip

\vskip 0.5 cm

\subheading{ Lemma 4.16}
{\it

$\overline{S}_{S\to T_S}^\parallel$ can be expressed by
$$
     \SHom_{\DMod_{T_S}}(\overline{S}_{S\to T_S}^\parallel, \CTS) =
     \{\ \psi\in \CTS\ | \ \Delta_{S \to T_S}^\parallel \psi =0 \},
$$
where 
$$
\split
\Delta_{S \to T_S}^\parallel &=  \Delta_S - \frac{1}{4} 
\frac{\tg_{T_S}^{\dalpha\dbeta}|_S }{g_S} \left(   \tr_{k\times 
k}(\gamma^{\alpha}_{\ \dalpha\beta}) \right) \left(\tr_{k\times 
k}(\gamma^{\alpha}_{\ \dbeta\beta})\right) \\
&+ \frac{1}{2}
\left( \frac{\tg_{T_S}^{\dalpha\dbeta}|_S }{g_{S}} 
        \tr_{k\times k}(\gamma^{\delta}_{\
         \dalpha\beta}\gamma^{\alpha}_{\ \dbeta\delta}) \right).
\endsplit
$$
}

\demo{Proof}
From the proposition 4.6, we have
$
\Delta_{S\to T_S} \equiv
 \Delta_{T_S}|_{q^\dalpha=0, \dalpha=k+1,\cdots,n}.
$
Proposition 3.19 gives
$$
 \split
    \Delta_{T_S} &= \tnablaT_\mu(\tg_{T_S})^{\mu\nu} 
    \tnablaT_\nu
    + \frac{1}{4} (\partial_\mu \log g_{T_S})( \partial_\mu 
    \tg_{T_S}^{\mu \nu })\\ &+ \frac{1}{4} \tg_{T_S}^{\mu\nu} ( 
    \partial_\mu \partial_\nu\log g_{T_S})+ \frac{1}{16} 
    \tg_{T_S}^{\mu\nu }
    (\partial_\mu \log g_{T_S})(\partial_\nu \log g_{T_S}).
\endsplit
$$
Since for $n, m \in \ZZ$, we have the relation
$
    \partial_\alpha^n \partial_\beta^m g_{T_S^\parallel}\equiv 
    \partial_\alpha^n \partial_\beta^m 
    g_S\quad \text{modulo } q^\dalpha
$,
 $\pi_{T_S}^*\Delta_{S}$ is given by,
$$
\split
    \pi_{T_S}^*\Delta_{S}&\equiv\tnablaT_\alpha 
   (\fg_{T_S})^{\alpha\beta} \tnablaT_\beta \\     
   &+ \frac{1}{4} (\partial_\alpha \log g_{T_S})( 
   \partial_\alpha g_{T_S}^{\alpha\beta}) 
   + \frac{1}{4} g_{T_S}^{\alpha\beta}
   ( \partial_\alpha \partial_\beta\log g_{T_S})\\
    &+ \frac{1}{16} 
   g_{T_S}^{\alpha\beta}(\partial_\alpha \log g_{T_S})
   (\partial_\beta \log
   g_{T_S})
\quad \text{modulo } q^\dalpha.
\endsplit
$$
By noting $\pi_{T_S}^*\Delta_{S}|_S\equiv \Delta_{S} $
and proposition 4.8 and by computing the remainder 
of $\Delta_{T_S}-\pi_{T_S}^*\Delta_{S}$, we obtain
$$
\split
\Delta_{S \to T_S} &=  \Delta_S + 
 \tnablaT_\dalpha(\tg_{T_S})^{\dalpha\dbeta} 
    \tnablaT_\dbeta |_S- \frac{1}{4} 
\frac{\tg_{T_S}^{\dalpha\dbeta}|_S }{g_S} \left(   \tr_{k\times 
k}(\gamma^{\alpha}_{\ \dalpha\beta}) \right) \left(\tr_{k\times 
k}(\gamma^{\alpha}_{\ \dbeta\beta})\right) \\
&+ \frac{1}{2}
\left( \frac{\tg_{T_S}^{\dalpha\dbeta}|_S }{g_{S}} 
        \tr_{k\times k}(\gamma^{\delta}_{\
         \dalpha\beta}\gamma^{\alpha}_{\ \dbeta\delta}) \right).
\endsplit
$$
From the proposition 4.6 again, we have
$
\Delta_{S \to T_S}^\parallel \equiv
 \Delta_{T_S}|_{ \tnablaT_\dalpha^\dalpha=0, \dalpha=k+1,\cdots,n}.
$ \qed
\enddemo

\vskip 0.5 cm

\subheading{Proof of Theorem 4.4}

Direct computations of the operator in 4.16 
leads us to obtain the theorem 4.4 noting the examples 4.14. \qed

\enddocument